\documentclass{elsarticle}
\usepackage[left=1in,top=1in,right=1in,bottom=1in]{geometry}
\usepackage{graphicx}
\usepackage{placeins}
\usepackage{subfig}
\usepackage{amsmath}
\usepackage{amssymb}
\usepackage{empheq}

\usepackage[colorlinks=true, allcolors=blue]{hyperref}
\usepackage[capitalise,sort&compress]{cleveref}
\biboptions{sort&compress}
\newtheorem{remark}{Remark}

\begin{document}

\begin{frontmatter}

\journal{arXiv}

\title{Pressure Stability in Explicitly Coupled Simulations of Poromechanics with Application to CO$_2$ Sequestration}

\author[1]{Ryan M. Aronson\corref{cor1}}
\ead{ryan.aronson@totalenergies.com}

\author[4]{Pavel Tomin}

\author[2]{Nicola Castelletto}

\author[1]{François P. Hamon}

\author[2]{Joshua A. White}

\author[5]{Hamdi A. Tchelepi}

\cortext[cor1]{Corresponding author}

\affiliation[1]{organization={TotalEnergies E\&P Research and Technology},
country={United States}}

\affiliation[4]{organization={Chevron Technical Center},
country={United States}}

\affiliation[2]{organization={Atmospheric, Earth, and Energy Division, Lawrence Livermore National Laboratory},
country={United States}}


\affiliation[5]{organization={Energy Science and Engineering, Stanford University},
country={United States}}

\begin{abstract}
We study in detail the pressure stabilizing effects of the non-iterated fixed-stress splitting in poromechanical problems which are nearly undrained and incompressible.
When applied in conjunction with a spatial discretization which does not satisfy the discrete inf-sup condition, namely a mixed piecewise linear - piecewise constant spatial discretization, the explicit fixed-stress scheme can have a pressure stabilizing effect in transient problems.
This effect disappears, however, upon time step refinement or the attainment of steady state.
The interpretation of the scheme as an Augmented Lagrangian method similar to Uzawa iteration for incompressible flow helps explain these results.
Moreover, due to the slowly evolving solution within undrained seal regions, we show that the explicit fixed-stress scheme requires very large time steps to reveal its pressure stabilizing effect in examples of geologic CO$_2$ sequestration.
We note that large time steps can result in large errors in drained regions, such as the aquifer or reservoir regions of these examples, and can prevent convergence of nonlinear solvers in the case of multiphase flows, which can make the explicit scheme an unreliable source of pressure stabilization.
We conclude by demonstrating that pressure jump stabilization is as effective in the explicit fixed-stress setting as in the fully implicit setting for undrained problems, while maintaining the stability and convergence of the fixed-stress split for drained problems. 
\end{abstract}

\begin{keyword}
Poromechanics \sep Explicit fixed-stress \sep Inf-sup stability \sep Pressure stabilization 
\end{keyword}

\end{frontmatter}

\section{Introduction}
As the world continues its progress towards green energy, geologic carbon sequestration will be critical in mitigating the environmental impact of industries that are more difficult or will take longer to decarbonize.
Predictive numerical simulation techniques are critical for the safety of these processes, ensuring that the subsurface pressure increases resulting from carbon dioxide injection do not reactivate existing faults or compromise the structural integrity of the storage formation itself.
Generating accurate numerical simulations is challenging, however, as solutions can involve a number of different, interacting physics such as multiphase fluid flow in porous media, geomechanics, thermal effects, and more. 

In this work we focus on single-phase and compositional poromechanics, or problems in which the deformation of a subsurface formation is coupled with the flow of one or more fluids within.
By compositional, we mean that mass transfer can occur between phases in the multiphase case.
Many numerical schemes used to solve these types of problems can be divided into three categories.
The first category are the monolithic, or fully implicit, algorithms in which the equations governing flow and mechanics are solved simultaneously.
These schemes can be made to perform well regardless of the strength of the coupling between the flow and mechanical processes, and can be proven to be unconditionally stable under appropriate assumptions \cite{garipov2018unified}.
However, fully implicit approaches can be highly complex and the solution of the resulting linear systems can require novel preconditioners specialized for multiphysics systems \cite{bui2020scalable, bui2021multigrid}.

The second broad category of algorithms we consider are the so-called sequential-implicit, or iterative-sequential, schemes.
Instead of advancing the mechanics and flow unknowns in time simultaneously, sequential implicit methods perform a splitting in order to solve a series of decoupled flow and mechanics problems.
These schemes are attractive, as they can repurpose existing, highly-optimized single-physics solvers to solve multiphysics problems.
Perhaps the most common sequential implicit method in reservoir engineering is the fixed-stress method \cite{settari1998coupled, settari2001advances, mainguy2002coupling, kim2011stability, kim2011stability_spe}, which has been shown to be convergent to the fully implicit solution provided sufficient iterations are performed \cite{settari1998coupled, kim2011stability, mikelic2013convergence, castelletto2015accuracy}.
The convergence rate can depend on the properties of underlying coupled system, however. 
For instance, the convergence of the sequential iterations can be slow in tightly coupled problems \cite{castelletto2015accuracy}, and various approaches have been proposed to accelerate convergence \cite{jeannin2007accelerating, storvik2018optimization, storvik2019optimization,waziri2024nonlinAccel}.

We refer to the final category of solvers as explicit, or non-iterative sequential schemes.
Like the previous set of algorithms, these approaches solve the flow and mechanics problems separately, and thus can make use of simpler techniques for each.
Unlike sequential implicit schemes, however, no iterations are performed within an individual time step and thus convergence to the fully implicit solution only occurs in the limit of mesh and time step refinement in general.
The fixed-stress scheme also functions as an explicit scheme when iterations are not performed (or when a fixed number of iterations are performed per time step) \cite{settari2001advances, kim2011stability}, and this will be the scheme considered in this work.
Finally, we remark that strategies in which only one of the subprocesses drives the other (and not vice-versa) could also be defined as explicitly coupled schemes.
It is common in the reservoir engineering community, for example, to ignore the coupling of the mechanics back onto the flow problem.
In this work, however, we will ignore these one-way coupled schemes and instead focus on the two-way coupled approaches mentioned above.

While the behavior of all of these types of coupling methods is perhaps well known in typical reservoir engineering contexts, application to CO$_2$ sequestration problems can require additional considerations.
In particular, formations which are attractive for carbon storage are usually bounded above by very low permeability caprocks, as these structures prevent the CO$_2$ plume from migrating upwards indefinitely due to buoyancy in a process known as structural trapping \cite{ipcc}.
It is well known that, when combined with solid grain and fluid components which are nearly incompressible (like water or brine), the governing equations of poromechanics reduce to a saddle-point, or constrained optimization, system similar to the equations governing incompressible flow \cite{zienkiewicz1990static}.
Accurate numerical simulation of these types of problems requires spatial discretizations which satisfy the Ladyzhenskaya–Babu\v{s}ka–Brezzi, or inf-sup, stability criterion \cite{brezzi2012mixed}.
The spatial discretization of choice within the reservoir simulation community combines a standard Finite Element (FE) method for mechanics with a Finite Volume (FV) approach for fluid flow and transport \cite{settari1998coupled,kim2011stability_spe,prevost2014,GarKarTch16}.
In particular, lowest-order piecewise-continuous interpolation is typically used for the displacement field, while flow quantities -- i.e., pressure and component densities -- are approximated in a piecewise-constant manner in each cell.
Such a choice does not satisfy the inf-sup criterion, and spurious pressure oscillations have been observed in simulations of CO$_2$ sequestration including caprocks using a fully implicit coupling strategy \cite{aronson2023pressure, aronson2024FS}.

The behavior of sequential methods (both implicit and explicit) in nearly undrained regimes is much less explored in the literature.
Upon first consideration, it seems that the inf-sup condition should not be required in these cases, as the saddle-point problem is never explicitly formed or solved, and this was the assertion in \cite{yoon2018spatial}.
However, later studies required the assumption of an inf-sup stable discretization when formulating a proof of convergence which was also valid in the undrained limit \cite{storvik2018optimization, storvik2019optimization}.
In our previous work \cite{aronson2024FS}, we clarified and reconciled these seemingly contradictory results by explaining that the key factor which determines if spurious pressure modes appear is the splitting error of the sequential scheme, not the invertibility of the discrete matrices appearing within the formulation.
This was inspired by similar debates and results in the incompressible flow community \cite{guermond2006overview, guermond1998stability}.
The focus of \cite{aronson2024FS} was a demonstration that, when sequential implicit schemes are iterated to convergence within a time step, the exact same pressure oscillations appeared as those appearing in a fully implicit solution.
The explicit case, where iterations are not included (and thus the splitting error is not driven to zero within time steps), was not studied in detail.
By considering similar splitting schemes developed for incompressible flow, we see that studies of the pressure stability of explicit coupling schemes must be carefully designed, as the stability properties can change under time step refinement, for example \cite{guermond2006overview, guermond1998stability, badia2007convergence}.

Such a study of the explicit fixed-stress method is the main goal of the current work.
After reviewing the governing equations of both single-phase and compositional poromechanics we detail the explicit fixed-stress method.
We then focus on the single-phase case when undrained conditions are approached.
The goal of this section is not to attempt a detailed proof of stability for discretizations which do not satisfy the inf-sup condition (which can be done in some cases for incompressible flow \cite{badia2007convergence}), but instead detail the operations performed in order to develop some intuitive expectations and motivate our numerical studies.
We show that the fixed-stress method in the undrained regime functions very similarly to Augmented Lagrangian methods \cite{fortin2000augmented}, in particular the classical Uzawa iteration for incompressible flow problems \cite{robichaud1990iterative}.
This provides an alternative, but complementary view to our previous work \cite{aronson2024FS}, where the role of splitting error was emphasized, and gives intuition on situations where the explicit fixed-stress method may be able to provide an implicit regularizing effect on the pressure.
We follow this with a series of numerical tests inspired by the analysis, where we find that, while the stabilizing effect does appear in certain cases, it is difficult to see in realistic cases of CO$_2$ sequestration due to specifics of their setup.
We then give a brief review of pressure jump stabilization \cite{hughes_1987, silvester1990stabilised, berger2015stabilized, camargo2021macroelement, aronson2023pressure}, which can be used to remove the remaining spurious oscillations, and was shown to be effective in the iterative fixed-stress setting \cite{aronson2024FS}.
We also include a brief discussion on how the inclusion of jump stabilization affects the analytical studies of the fixed-stress method performed in \cite{kim2011stability} for drained problems, highlighting that stability in time and convergence of the explicit scheme are preserved for single-phase cases in one spatial dimension.
Finally, we return to our examples to demonstrate numerically the effectiveness of the stabilization in the explicit fixed-stress setting.

\section{Model Problems}

In this section we review the governing equations for isothermal compositional multiphase flow and transport in a deformable porous medium assuming an arbitrary number of components $n_c$ and phases $n_p$.
Note that a thorough derivation of the model is beyond the scope of the present work -- additional details can be found for example in \cite{biot1941general, coussy2004poromechanics, terzaghi}.
As our application of interest is in geologic carbon sequestration, in our numerical results we will focus on the two phase ($n_p = 2$), two component CO$_2$-brine system ($n_c = 2$) where the CO$_2$ can be present in either phase while the brine component is limited to the liquid phase.
We will also specialize the general formulation for the single-phase ($n_p = 1$), single-component ($n_c = 1$) case that will mainly be used for analysis and developing understanding of the behavior of the fixed-stress method.

We consider a displacement ($\mathbf{u}$), pore pressure ($p$), global component densities ($\boldsymbol{\rho}_c$) formulation.
Here, $\boldsymbol{\rho}_c = \{ \rho_c \}_{c  = 1, \ldots, n_c}$ is the vector that collects the global density of each mobile component $\rho_c = \sum_{\ell=1}^{n_p} x_{c\ell} \rho_{\ell} s_{\ell} $, with $x_{c\ell}$, $\rho_{\ell}$ and $s_{\ell}$ the phase component fraction, density and saturation, respectively.
In this work we ignore any capillarity effects, thus $p$ at a point is the same across all phases.
Denoting by $\Omega \subset \mathbb{R}^3$ and $\mathcal{I} = (0, t_{\text{final}}] $ the domain and the time interval of interest, the initial-boundary value problem of interest read as follows \cite{coussy2004poromechanics}:
Find 
$\mathbf{u}: \Omega \times \mathcal{I} \mapsto \mathbb{R}^3$,
$p: \Omega \times \mathcal{I} \mapsto \mathbb{R}$,
$\boldsymbol{\rho}_c: \Omega \times \mathcal{I} \mapsto \mathbb{R}^{n_c}$
such that:
\begin{subequations}
\begin{align}
  & - \nabla \cdot \boldsymbol{\sigma} ( \mathbf{u}, p )
  = \mathbf{f} ( \mathbf{u}, p, \boldsymbol{\rho}_c )
  &
  &
  \text{(quasi-static linear momentum balance)},
  \label{eq:biot_multiphase_momentum}\\
  & \dot{m}_c( \mathbf{u}, p, \rho_c )
  + \nabla \cdot \mathbf{j}_{m,c}( p, \boldsymbol{\rho}_c )
  = q_{m,c}( p, \boldsymbol{\rho}_c )
  &
  &
  \text{(mass balance for component $c$, with $c \in \{1, \ldots, n_c \}$)},
  \label{eq:biot_multiphase_mass}\\
  & \sum_{\ell = 1}^{n_p} s_{\ell}( p, \boldsymbol{\rho}_c )
  = 1
  &
  &
  \text{(saturation constraint)},
  \label{eq:biot_multiphase_saturation_constraint}
\end{align}
\label{eq:biot_multiphase}\null
\end{subequations}
subject to suitable boundary and initial conditions.
In \cref{eq:biot_multiphase}, the following symbols, variables, and constitutive relationships are introduced:
\begin{itemize}

  \item $\boldsymbol{\sigma} ( \mathbf{u}, p ) = ( \boldsymbol{\sigma}^{\prime}  - \mathbf{b} p )$ is the total stress tensor, with $\boldsymbol{\sigma}^{\prime} = ( \mathbb{C}_{dr} : \boldsymbol{\varepsilon} ( \mathbf{u} )$ the effective stress tensor, $\mathbb{C}_{dr}$ the fourth-order stress–strain tangent tensor, $\boldsymbol{\varepsilon} ( \mathbf{u} ) = \nabla^s \mathbf{u}$ the linearized second-order strain tensor, $\nabla^s$ the symmetric gradient operator, and $\mathbf{b}$ is the Biot coefficient tensor. As the saddle-point structure we are interested in investigating is not associated with any particular type of constitutive response of the rock matrix, we consider only linear elastic isotropic behavior in this work. Hence:
  \begin{align}
    \mathbb{C}_{dr}
    &
    = K_{dr} (\mathbf{1} \otimes \mathbf{1} )
    + 2G \left( \mathbb{I} - \frac{1}{3} \mathbf{1} \otimes \mathbf{1}  \right),
    &
    \mathbf{b}
    &
    = b \mathbf{1}, 
    &
    b 
    &= 
    1 - \frac{K_{dr}}{K_s},
  \end{align}
with $K_{dr}$ and $K_s$ the bulk modulus of the drained skeleton and the solid phase, respectively, $G$ the shear modulus, and $\mathbf{1}$ and $\mathbb{I}$ the second- and fourth-order identity tensor, respectively. We also define the following auxiliary terms: (i) the volumetric total stress $\sigma_v = \frac{1}{3}\text{tr}(\boldsymbol \sigma)$, (ii) the volumetric effective stress $\sigma'_v = \frac{1}{3}\text{tr}(\boldsymbol \sigma^{\prime})$, and (iii) the volumetric strain $\epsilon_v = \text{tr}(\boldsymbol \varepsilon)$, where $\text{tr}(\bullet)$ denotes the trace of quantity $(\bullet)$. Note that drained bulk modulus of the solid skeleton relates volumetric stress and strain via $\sigma'_v = K_{dr}\epsilon_v$;

  \item $\mathbf{f} ( \mathbf{u}, p, \boldsymbol{\rho}_c ) = \rho( \mathbf{u}, p, \boldsymbol{\rho}_c ) \mathbf{g}$ denotes the body force vector, where $\rho = (1 - \phi( \mathbf{u}, p )) \rho_s + \phi( \mathbf{u}, p ) \rho_T$ is the mixture density and $\mathbf{g}$ is the gravity vector, with $\rho_s$ the solid phase density, $\rho_T = \sum_{c = 1}^{n_c} \rho_c$ the total component mass density, and $\phi$ the porosity. The porosity change from a reference porosity $\phi_0$ is modeled in incremental form using the relationship \cite{coussy2004poromechanics}:
\begin{align}
  \phi( \mathbf{u}, p ) &= \phi_0 + b ( \epsilon_v( \mathbf{u} ) - \epsilon_{v,0} ) + \frac{1}{N} ( p - p_0 ),
  &
  \frac{1}{N} &= \frac{(b-\phi_0)}{K_s};
\end{align}

  \item $m_c( \mathbf{u}, p, \rho_c ) = \phi ( \mathbf{u}, p ) \rho_c$ denotes the mass per unit volume for component $c$. In \cref{eq:biot_multiphase_mass}, the superposed dot, $\dot{(\bullet)}$, indicates the derivative with respect to time of quantity $(\bullet)$;

  \item $\mathbf{j}_{m,c} (p, \boldsymbol{\rho}_c  )$ is the mass flux for component $c$, which accounts for all the phases flux contributions based on a traditional multiphase flow extension of Darcy's law, namely
  \begin{align}
  & \mathbf{j}_{m,c} (p, \boldsymbol{\rho}_c )
  = - \sum_{\ell = 1}^{n_p}
  x_{c\ell}(p, \boldsymbol{\rho}_c )
  \rho_{\ell}(p, \boldsymbol{\rho}_c )
  \frac{k_{r\ell}(p, \boldsymbol{\rho}_c )}{\mu_{\ell}(p, \boldsymbol{\rho}_c )}
  \boldsymbol{\kappa} \cdot ( \nabla p - \rho_{\ell}(p, \boldsymbol{\rho}_c ) \mathbf{g} ),
  \label{eq:darcy_multiphase}
  \end{align}
   with $k_{r\ell}$ and $\mu_{\ell}$ the relative permeability and viscosity of phase $\ell$, respectively, $\boldsymbol{\kappa}$ the second-order permeability tensor, and $\nabla$ the gradient operator.
  The fluid properties are updated using the closure relations
  \begin{align}
    \sum_{\ell=1}^{n_p} s_{\ell} &= 1,
    &
    \sum_{c=1}^{n_c} x_{c\ell} &= 1,  \qquad \ell = 1, \ldots, n_p,
  \end{align}
  and thermodynamic equilibrium constraints expressed in terms of the fugacity $f_{c,\ell}(p, T, x_{c,\ell})$ of component $c$ in phase $\ell$ :
  \begin{equation}
    f_{c\ell}(p, T, x_{c\ell}) - f_{ck}(p, T, x_{ck}) = 0 \qquad \forall \ell \neq k, \quad c = 1, ... , n_c,
  \end{equation}
  with $T$ representing temperature.
  Assuming a constant known temperature, the update procedure involves the following steps:  
  First, the phase volume fractions ($s_{\ell}$) and phase component fractions ($x_{c\ell}$) are computed as a function of pressure ($p$), and total component fractions ($z_c = \rho_c / \rho_T$).
  Then, phase densities ($\rho_{\ell}$) and phase viscosities ($\mu_{\ell}$) are computed as a function of pressure and the updated phase component fractions, and phase relative permeabilities ($k_{r\ell}$) are computed as a function of the updated phase volume fractions ($s_{\ell}$).
  More detail can be found in the online documentation of GEOS \cite{GEOS}, which we use to conduct the numerical studies in this work.

  \item $q_{m,c}(p, \boldsymbol{\rho}_c  ) = \sum_{\ell = 1}^{n_p} x_{c\ell} \rho_{\ell} q_{\ell}$ denotes the mass source/sink per unit volume for component $c$, with $q_{\ell}$ the volumetric flow rate for phase $\ell$. In this work, simple, element-wise constant source fluxes of CO$_2$ are used to model injection, but terms can also be introduced to model wells based on inflow-performance relationships that depend on $p$ and $\boldsymbol{\rho}_c$ \cite{lie19}.

\end{itemize}

For single-phase, single-component poromechanics, the initial-boundary problem \cref{eq:biot_multiphase} only requires \cref{eq:biot_multiphase_momentum,eq:biot_multiphase_mass}, with the local saturation constraint becoming irrelevant.
The primary unknowns reduce to $\mathbf{u}$ and $p$, with $\rho_f(p) = \rho_c = \rho_T = \rho_{\ell}$ and $\mu_f(p) = \mu_{\ell}$ the density and the viscosity of the single-phase, single-component fluid that depend on $p$ based on suitable constitutive equations ($c =\ell=1$).
Note that both the component fraction $x_{c\ell}$ and the relative permeability coefficient are equal to 1.
Following \cite{kim2011stability}, in which different formulations of the mass equations are used in the development of sequential implicit schemes, we also consider this alternative statement of \cref{eq:biot_multiphase_mass} assuming a slightly compressible single-phase, single-component fluid
\begin{equation}
    \left( \frac{1}{M} + \frac{b^2}{K_{dr}} \right) \dot{p} + \frac{b}{K_{dr}}
    \dot{\sigma_v} + \nabla \cdot \mathbf{v} =  q,
    \label{eq:single_phase_mass_FS}
\end{equation}
where $M = \frac{1}{N} + \frac{\phi}{K_f}$ denotes the Biot modulus, $K_f$ the fluid bulk modulus, and $\mathbf{v}$ the single-phase Darcy flux.

\subsection{Undrained Flows with Incompressible Materials}

While the above sets of equations can be used to model the poromechanical response of many porous materials, it is important to note how their character changes when both the solid grains and fluid are incompressible and the permeability magnitude is equal to zero.
These conditions can occur, for instance, in the caprock structures providing structural trapping in cases of geologic CO$_2$ sequestration. 

For the single-phase case, starting with \cref{eq:biot_multiphase_mass} and applying the assumptions stated above (including that one cannot inject into or produce from a zero permeability formation) leads to
\begin{equation}
     \rho_f \dot{ \phi } = 0.
\end{equation}
Noting that the assumptions also imply $b = 1$, this further reduces to 
\begin{equation}
      \dot{ \epsilon_v} = \nabla \cdot \dot{\mathbf{u}} = 0.
      \label{eq:incompress_constraint}
\end{equation}
When combined with \cref{eq:biot_multiphase_momentum} describing mechanical force equilibrium, we see that the poromechanical problem has a saddle-point structure, equivalent to the governing equations of incompressible Stokes flow:
\begin{subequations}
  \begin{empheq}[left=\empheqlbrace]{align}
     \nabla \cdot (\mathbb{C}_{dr} : (\nabla^s \mathbf{u})) - b \nabla p &  = 0,\\
     \nabla \cdot \dot{\mathbf{u}} &  =  0.
  \end{empheq}
  \label{eq:continuous_saddle}
\end{subequations}
In particular, the mass equation no longer takes the form of an evolution equation for the pressure, but instead takes the form of a constraint equation, with the pore pressure acting as a Lagrange multiplier enforcing incompressibility of the solid deformations.
Also note that, if the initial displacement solution is divergence free (as would be the case for a zero initial condition, for example), then the constraint equation can also be written as $\nabla \cdot \mathbf{u} = 0$, meaning the displacement is divergence free instead of its velocity.
\begin{remark}
    For brevity we do not include any explicit discussion of the multiphase case, but we note that the reduction to a saddle-point problem in the case of an undrained, incompressible problem is very similar; conclusions regarding the stability of the fixed-stress method in the single-phase, undrained, incompressible case also apply to multiphase scenarios. We refer to \cite{camargo2021macroelement} for more detail on the immiscible multiphase setting and to \cite{aronson2023pressure, aronson2024FS} for the compositional setting. 
\end{remark}

\section{Spatial Discretization}
\label{sec:spatial_discretization}

To discretize the above governing equations in space, we elect to use a mixed finite element - finite volume formulation where the displacement degrees of freedom are discretized with lowest-order piecewise-continuous nodal finite elements and the remaining degrees of freedom are discretized using piecewise constant, cell-centered values, where the numerical flux is computed using a standard Two Point Flux Approximation (TPFA) \cite{aziz1979petroleum}.
For brevity, we defer a detailed discussion to our previous work \cite{aronson2024FS}, which details the spatially discrete versions of the equations which are also valid in this work.

Here we will simply state that, for the single-phase case, the spatially discrete equations can be written as a block system of differential-algebraic equations of the form
\begin{equation}
    \begin{bmatrix}
    A & -B^T \\
      & T
    \end{bmatrix}
    \begin{bmatrix}
    u \\
    p
    \end{bmatrix}
    +
    \begin{bmatrix}
      &  \\
    B & M
    \end{bmatrix}
    \begin{bmatrix}
    \dot{u} \\
    \dot{p}
    \end{bmatrix}
    =
    \begin{bmatrix}
    Q_u \\
    Q_p
    \end{bmatrix}.
    \label{eq:block_system_DAE}
\end{equation}
Here $A$ represents the elastic stiffness matrix, $B$ and $B^T$ are the discrete divergence and gradient operators, respectively, $T$ contains the contributions from the flux terms, $M$ contains contribution from the accumulation term, and $Q_u$ and $Q_p$ represent any forcing terms.
Time discretization is then done with a backward Euler method, either applied to the entire coupled problem shown in \cref{eq:block_system_DAE} in the case of a fully implicit method, or to each single physics problem independently in the case of a fixed-stress solver, which is the topic of \cref{sec:explicit_FS_split}. 

\begin{remark}
    The compositional multiphase case results in a nonlinear system of equations which is solved with a Newton-Raphson approach, but the Jacobian at any instant in time is structured similarly to Eq. \cref{eq:block_system_DAE}, but includes more unknowns. In our implementation, we utilize pressure and global component densities as primary unknowns, though other choices are possible \cite{voskov2012comparison}. Secondary unknowns can be eliminated via static condensation as in \cite{bui2021multigrid}. 
\end{remark}

\subsection{Undrained Flows with Incompressible Materials}

Similar to the fully continuous setting, we also highlight the changes that occur when the permeability is small and the fluid and solid grains are incompressible.
The discrete analog of \cref{eq:continuous_saddle} can be written as
\begin{equation}
    \begin{bmatrix}
    A & -B^T \\
    B & 0
    \end{bmatrix}
    \begin{bmatrix}
    u \\
    p
    \end{bmatrix}
    =
    \begin{bmatrix}
    Q_u \\
    Q_p
    \end{bmatrix}.
    \label{eq:discrete_saddle}
\end{equation}
The spatial discretizations used to solve \cref{eq:discrete_saddle} in a fully implicit manner must be carefully selected for stability, in particular satisfying the discrete inf-sup condition \cite{brezzi2012mixed}.
Results obtained with unstable element pairs will exhibit spurious oscillations in the discrete pressure field sometimes referred to as checkerboard oscillations \cite{sani_causesandcure}.
A spurious pressure mode is one which is not identically zero, $p \neq 0$, but satisfies $B^T p = 0$, meaning the discrete gradient operator has a nontrivial null space. 
The remainder of this work is devoted to determining if these modes also appear when utilizing the explicit fixed-stress method.

\begin{remark}
    It has been shown in the poromechanical context that these spurious modes will appear even if incompressibility and zero permeability conditions are merely approached and not attained exactly, in both single and multiphase poromechanics problems.
    Moreover, realistic values used in simulations of CO$_2$ sequestration satisfy these conditions closely enough for pressure oscillations to appear \cite{aronson2023pressure, aronson2024FS}.
\end{remark}

\section{The Explicit Fixed-Stress Split}
\label{sec:explicit_FS_split}

With the governing equations of both single and multiphase poromechanics described, we now consider discretization in time using the fixed-stress splitting. We focus on the explicit case where no iterations are taken within time steps, so that each time step consists of one flow solve and one mechanics solve.

The explicit fixed-stress split first solves the flow problem with the rate of the mean total stress held constant \cite{kim2011stability}.
More formally, when advancing the flow problem from time $n$ to time $n+1$, it is assumed that  
\begin{equation}
     \sigma^{n+1}_{v} - \sigma^{n}_{v} = \sigma^{n}_{v} - \sigma^{n-1}_{v}.
     \label{eq:FS_split}
\end{equation}
In the compositional setting, the flow problem is nonlinear and its solution requires an iterative Newton-Raphson solver, while the problem in the single-phase case is linear.
Once the the flow solution is known at time $n+1$, it is used to determine the final displacement at time $n+1$.
As we restrict ourselves to linear elasticity in this work, the mechanics problem is always linear and does not require multiple Newton iterations to solve.

The fixed-stress splitting has been shown to be convergent to the fully implicit solution when iterated to convergence \cite{kim2011stability, castelletto2015accuracy, settari1998coupled}, and more importantly for this work, the explicit scheme (where no iterations are taken) also converges under refinement \cite{kim2011stability, mikelic2013convergence}.
The fixed-stress split can also be written in the form of a modification to the porosity field, as discussed in \cite{kim2011stability_spe, garipov2018unified}.
In this case the porosity is updated according to 
\begin{equation}
    \phi^{n+1} = \phi^0 + \frac{(b-\phi^0)(1-b)}{K_{dr}} (p^{n+1} - p_0) + b\epsilon_v^n + \frac{b^2}{K_{dr}}(p^{n+1} - p^n),
\end{equation}
and this value is passed between the solvers.
This view of the fixed-stress split reveals that it is equivalent to the decoupling strategies developed previously in \cite{settari1998coupled, settari2001advances}.
Moreover, this approach is easily adapted for both single-phase and compositional poromechanics, and so it is the approach we utilize in our numerical implementation.

For the purposes of understanding its behavior when applied to the undrained saddle-point problem, however, we will use the explicit expressions of the equations for the single-phase poromechanical case developed in \cite{kim2011stability}.
We will see through numerical tests that these conclusions still hold in the compositional setting, and in particular in realistic simulations of CO$_2$ injection.
The starting point of the scheme is \cref{eq:single_phase_mass_FS}.
Ignoring forcing terms, using Darcy's law (\cref{eq:darcy_multiphase}) to replace the velocity term with pressure, discretizing in time with a backward Euler scheme with time step of size $\delta t$, and applying the fixed-stress assumption (\cref{eq:FS_split}) yields
\begin{equation}
    \left(\frac{1}{M} + \frac{b^2}{K_{dr}} \right) \frac{p^{n+1} - p^{n}}{\delta t} + \left( \frac{b}{K_{dr}}\right) \frac{\sigma_v^{n} - \sigma_v^{n-1}}{\delta t} - \nabla \cdot \left( \frac{\mathbf{k}}{\mu} \nabla p^{n+1} \right) = 0. 
    \label{eq:FS_mass}
\end{equation}
\cref{eq:FS_mass} is used to solve for the pressure $p^{n+1}$ at the new discrete time level $t^{n+1}$.
Ignoring external forcing again for simplicity, the displacement $u^{n+1}$ is then determined using
\begin{equation}
    \nabla \cdot (\mathbb{C}_{dr} : (\nabla^s \mathbf{u}^{n+1})) - b \nabla p^{n+1} = 0,
    \label{eq:FS_mech}
\end{equation}
and then the process is repeated for the subsequent time step.

\subsection{Application to Saddle-Point Problems}


Of interest to us in this work is the behavior of the non-iterated fixed-stress splitting as the governing equations approach a saddle-point system. In this case, standard fractional step or splitting theory is not applicable \cite{janenko1971method}.
We focus on the spatially continuous setting for simplicity, but we note that we will discretize in space according to \cref{sec:spatial_discretization}.
With the assumptions of perfectly incompressible solid grains and fluid and exactly zero permeability, the fixed-stress mass equation (\cref{eq:FS_mass}) becomes
\begin{equation}
    \left(\frac{1}{K_{dr}} \right) (p^{n+1} - p^{n}) + \left( \frac{1}{K_{dr}}\right) (\sigma_v^{n} - \sigma_v^{n-1}) =0,
\end{equation}
which can be simplified to 
\begin{equation}
    (p^{n+1} - p^{n}) = (p^{n} - p^{n-1}) - K_{dr}\nabla \cdot (\mathbf{u}^n - \mathbf{u}^{n-1}).
    \label{eq:FS_mass_undrained}
\end{equation}
Note the absence of the time step size $\delta t$ in \cref{eq:FS_mass_undrained}, which reflects the fact that the mass equation is no longer an evolution equation for pressure, but instead a constraint equation approximating \cref{eq:incompress_constraint} at the time-discrete level.

\subsection{Interpretation as an Augmented Lagrangian Method}

The easiest way to understand the behavior of the fixed-stress method in the undrained, incompressible setting is to simply walk through the algorithm at the time-discrete level:
\begin{enumerate}
\item At the first flow solve, the fixed-stress assumption reduces to $\sigma_v^1 = \sigma_v^0$, and so the resulting pressure field is simply given by
\begin{equation}
  p^1 = p^0.
\end{equation}
Then, assuming a non-constant mechanical forcing term $f$, the new displacement is computed via
\begin{equation}
  \nabla \cdot (\mathbb{C}_{dr} : (\nabla^s \mathbf{u}^{1}))  = f^1 + \nabla p^{1}.
\end{equation}
\item Continuing to the second time step, the fixed stress assumption is stated as $\sigma_v^2 - \sigma_v^1 = \sigma_v^1 - \sigma_v^0$, and so the flow problem is given by
\begin{equation}
  (p^{2} - p^{1}) = (p^{1} - p^0) - K_{dr} \nabla \cdot (\mathbf{u}^1 - \mathbf{u}^0).
\end{equation}
Using the result of the first flow solve, $p^1 = p^0$, then this simplifies to 
\begin{equation}
  (p^{2} - p^{1}) = - K_{dr} \nabla \cdot (\mathbf{u}^1 - \mathbf{u}^0),
\end{equation}
and the following mechanics step is given by 
\begin{equation}
  \nabla \cdot (\mathbb{C}_{dr} : (\nabla^s \mathbf{u}^{2}))  = f^2 + \nabla p^{2}.
\end{equation}
\item In the third time step, the new pressure is computed via
\begin{equation}
  (p^{3} - p^{2}) = (p^{2} - p^1) - K_{dr} \nabla \cdot (\mathbf{u}^2 - \mathbf{u}^1),
\end{equation}
which, by noting that $(p^{2} - p^{1}) = - K_{dr} \nabla \cdot (\mathbf{u}^1 - \mathbf{u}^0)$, becomes
\begin{equation}
    (p^{3} - p^{2}) = K_{dr} \nabla \cdot (\mathbf{u}^2 - \mathbf{u}^0).
\end{equation}
The mechanics step to be performed is then
\begin{equation}
    \nabla \cdot (\mathbb{C}_{dr} : (\nabla^s \mathbf{u}^{3}))  = f^3 + \nabla p^{3}.
\end{equation}
\end{enumerate}
Following this pattern, the fixed-stress flow solve for a general future time step can be written as
\begin{equation}
  p^{n+1} = p^n + K_{dr} \nabla \cdot (\mathbf{u}^n - \mathbf{u}^0),
  \label{eq:FS_penalty_form}
\end{equation}
meaning that the pressure is updated by the net volumetric displacement of the previous mechanics solution before being used as a forcing in the next mechanics solve.
This formulation resembles a penalty formulation for incompressible elasticity with the penalty parameter specified as the drained bulk modulus $K_{dr}$.
One can also interpret this as an attempt to define the new pressure to be the one which would make the previous displacement increment divergence-free.

In fact, if we assume the initial displacement to be divergence-free, we can re-write the coupled equations in a form which strongly resembles Augmented Lagrangian methods for incompressible flow \cite{fortin2000augmented}:
\begin{subequations}
  \begin{empheq}[left=\empheqlbrace]{align}
     \nabla \cdot (\mathbb{C}_{dr} : (\nabla^s \mathbf{u}^{n+1})) & = f^{n+1} + \nabla p^{n} + K_{dr}\nabla (\nabla \cdot \mathbf{u}^n)\\
     p^{n+1}& = p^n + K_{dr} \nabla \cdot (\mathbf{u}^n).
  \end{empheq}
  \label{eq:FS_uzawa}
\end{subequations}
The similarity of the fixed-stress split with Uzawa iteration \cite{robichaud1990iterative} (which is an Augmented Lagrangian approach) was noted in a remark in \cite{yoon2018spatial}, though they did not discuss its implication for pressure stability in detail.
Moreover, the authors did not note that the Uzawa method produces non-singular matrices at every stage, but still requires inf-sup stable discretizations to avoid pressure oscillations in the solution. This actually contradicts their main assertion about the stability of the fixed-stress method.
In fact, the fixed-stress splitting functions almost identically to an inexact Uzawa iteration \cite{bramble_inexactUzawa}.
Inexact Uzawa methods approximate the pressure Schur complement (usually with a scaled pressure mass matrix) so that no matrix inversion is needed, and this is analogous to the $b^2/K_{dr}$ diagonal term introduced in the mean fixed-stress split \cite{castelletto2015accuracy}.

The interpretation of the fixed-stress splitting as an Augmented Lagrange approach provides much clarity for the pressure stability of the method.
In a transient problem, which would be reflected in the above algorithm as time-varying forces $f^n$ (but could also be reflected with mass sources and sinks), the explicit fixed-stress split is the same as an Uzawa method with only one iteration.
However, Uzawa methods rely on the presence of multiple iterations within a time step to accurately enforce the incompressibility constraint.
When iterations are not taken, the incompressibility constraint is only accurately enforced in the limit of infinitely large penalty parameters, but this is fixed as $K_{dr}$ in the fixed-stress method.
Moreover, since the penalty term in \cref{eq:FS_uzawa} is based on a lagged displacement solution, the enforcement of the incompressibility constraint is further inaccurate in the explicit case if the displacement is rapidly changing.
This may be the case with large time steps, for example. Given this inaccuracy in the enforcement of the incompressibility constraint (which is equivalent to the splitting error discussed in \cite{aronson2024FS}), it may be reasonable to assume that pressure oscillations will not appear, as was concluded in \cite{yoon2018spatial}, albeit for the wrong reasons.

This conclusion is only true in truly transient problems, however.
In particular, consider using the above algorithm when the physical solution has reached a steady state.
If we continue to time step after the steady state solution has been reached, we should obtain solutions which are closer and closer to incompressible, as repeated applications of the penalty term will result in improved constraint satisfaction.
In this case the explicit scheme is functioning identically to the iterative scheme, and so we would expect pressure oscillations to return, given the results of \cite{aronson2024FS}.
We will verify this with our numerical results, and we note that this is common in industry simulations of CO$_2$ sequestration, where after injection the simulation is continued for many years to model the migration of the CO$_2$ plume.
We also note that, even in cases which are not fully steady, if the solution changes very slowly (when compared to the time step size), we may also expect to see pressure oscillations for the same reasons as above. In fact, we shall see that, depending on the time step size, pressure solutions in the burden regions of CO$_2$ sequestration problems seem to not change rapidly enough to avoid the spurious pressure modes. 

It is important to note that the above interpretation does not provide any rigorous proof of pressure stability (or lack thereof) in the explicit fixed-stress splitting with our choice of spatial discretization. Instead, we use the interpretations developed above as motivation in the design of the numerical studies in the following section. 

\begin{remark}
    The interpretation of the fixed-stress splitting as an Augmented Lagrange approach also agrees very well with the interpretations developed in \cite{aronson2024FS} for the iterative fixed-stress scheme.
    It is well-known that Augmented Lagrangian approaches provide no implicit pressure stabilizing effect, and thus if a non-inf-sup stable spatial discretization is applied, spurious pressure oscillations can still appear \cite{robichaud1990iterative, benzi_ALMPrecond}. 
    This is exactly what we found in our previous work \cite{aronson2024FS}, though there we used a motivation based on projection methods for incompressible flow and numerical splitting errors.
    Moreover, the large number of iterations required when solving undrained problems with non-inf-sup stable discretizations seen for the fixed-stress method in~\cite{aronson2024FS, storvik2018optimization}, was also seen when using Uzawa iteration to solve incompressible flow problems \cite{robichaud1990iterative}. 

    In~fact, as iterations are performed within a time step (and the constraint is enforced more accurately), the convergence rate of the Uzawa algorithm slows in simulations of incompressible flow with unstable spatial discretizations \cite{robichaud1990iterative}.
    We~also saw this change in convergence rate in our previous work \cite{aronson2024FS}.
\end{remark}

\section{Numerical Studies of Pressure Stability}

To gain a full understanding of the pressure stability properties of the fixed-stress method, we study three cases of increasing complexity and applicability to CO$_2$ sequestration.
We start by considering single-phase poromechanics in a cantilevered plate used in previous studies of pressure stability in mixed discretizations.
In this case we shall see cases where the fixed-stress method can exhibit improved pressure stability compared to the fully implicit method.
We then consider single-phase poromechanics in a staircase geometry to begin to understand stability properties in a case where only a burden region is undrained.
Then we conclude with a realistic CO$_2$ injection case, namely a compositional poromechanics simulation of injection into the High Island 24L formation.
This result demonstrates that the motivations developed using the single-phase case are also valid in more complex settings, and we highlight additional considerations that must be made for nonlinear flow problems.

We utilize GEOS \cite{GEOS, gross2024geos} to perform all of the numerical tests in this work. GEOS is an open-source, high-performance simulator with the capability to solve coupled poromechanical problems using both fully implicit as well as sequential (implicit and explicit) methodologies, which allows us to easily evaluate the pressure stability of the explicit fixed-stress split compared to the corresponding fully implicit method. As mentioned above, all nonlinear systems are solved via a Newton-Raphson approach, and linear systems are solved using GMRES with algebraic multigrid and multigrid reduction preconditioners \cite{bui2021multigrid}.

\subsection{Cantilever Plate}

\begin{figure}
\centering
\includegraphics[width=0.4\textwidth]{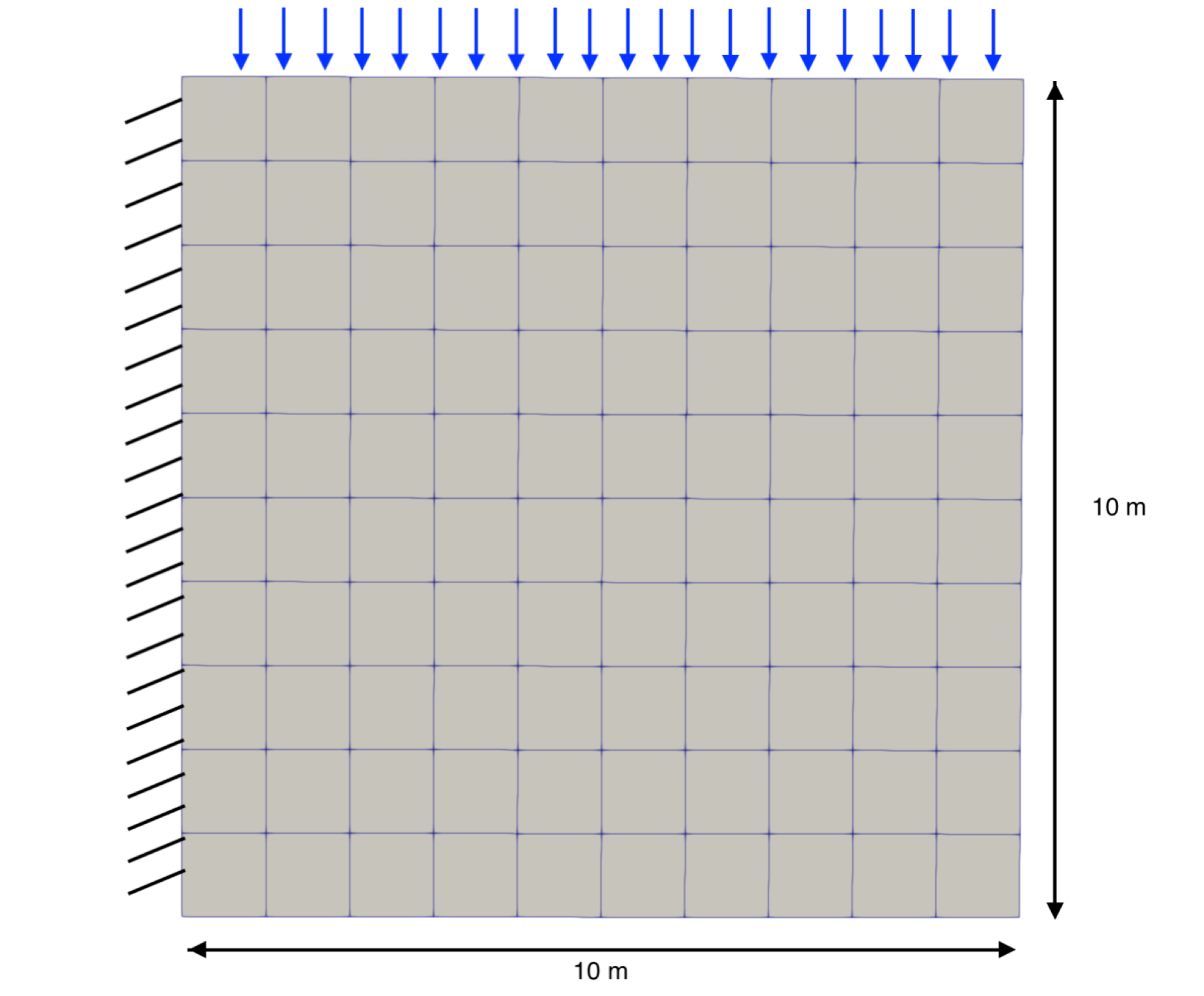}
\caption{Cantilever problem setup}
\label{fig:cant_setup}
\end{figure}

Our first example is inspired by the studies in \cite{phillips2009overcoming, borio2021hybrid, frigo2021efficient}, where it was used to assess the effectiveness of pressure stabilization strategies in other discretizations.
We consider a planar poroelastic solid given by the mesh in \cref{fig:cant_setup}.
The left edge of the skeleton is fixed and a sinusoidal loading is applied to the top with maximum magnitude of 100~N and a period of 10~days.
All other edges are free. The skeleton is assumed isotropic and linear elastic with drained bulk modulus equal to 5~GPa and Poisson's ratio $\nu = \frac{3K_{dr} - 2G}{2(3K_{dr} + G)} = 0.25$, and the solid grains are assumed to be exactly incompressible.

We consider single-phase flow within the skeleton, with no-flow boundary conditions applied on all four boundaries.
The permeability of the skeleton is defined to be isotropic and equal to zero everywhere, and the porosity is set to 0.05 everywhere.
The fluid is defined to have density 1000 kg/m$^3$ and zero compressibility.

\begin{figure}
\centering
\subfloat[$t = 3$ days]{\label{sfig:cant_fim_3}\includegraphics[width=.4\textwidth]{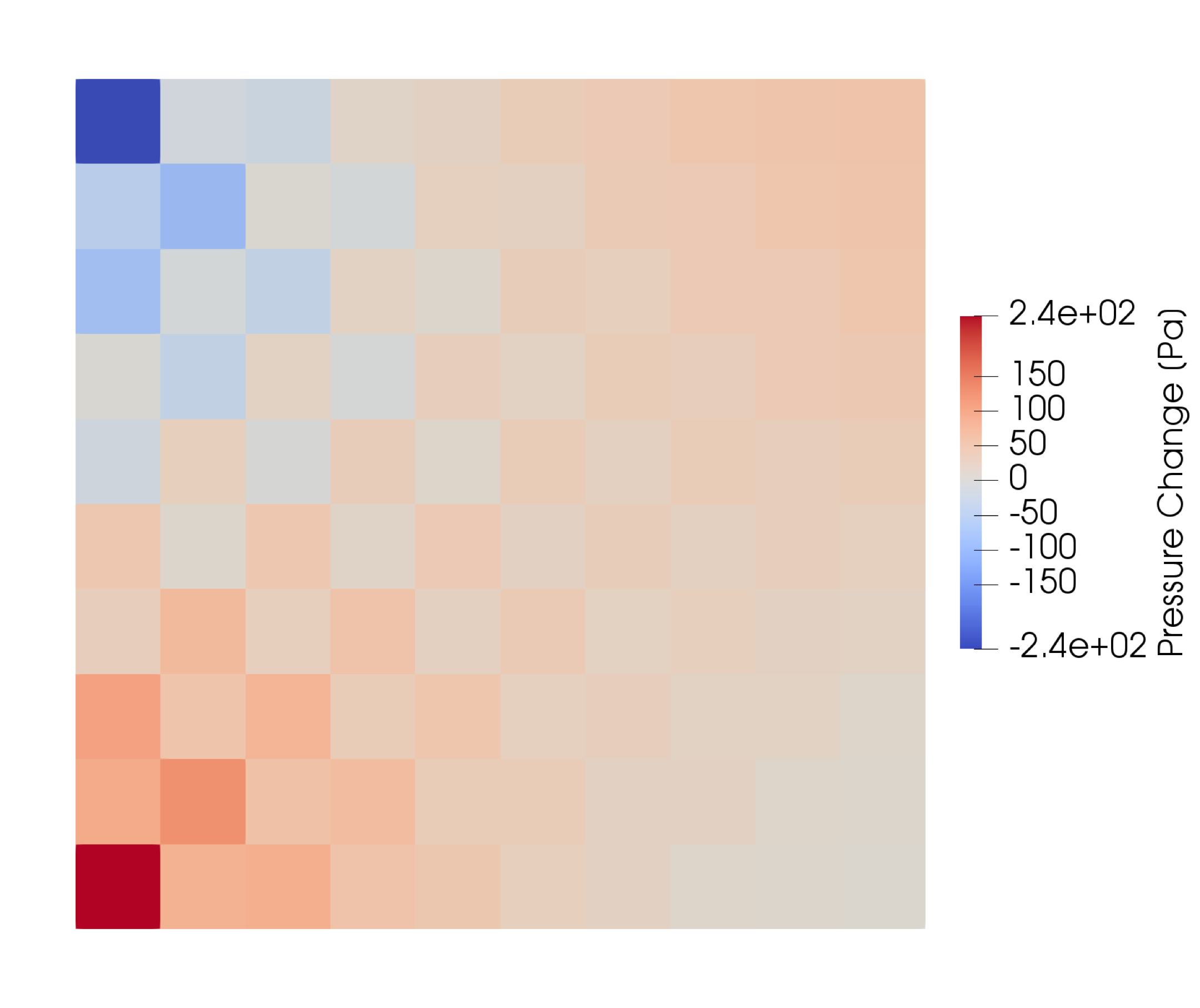}}\hfill
\subfloat[$t = 7$ days]{\label{sfig:cant_fim_7}\includegraphics[width=.4\textwidth]{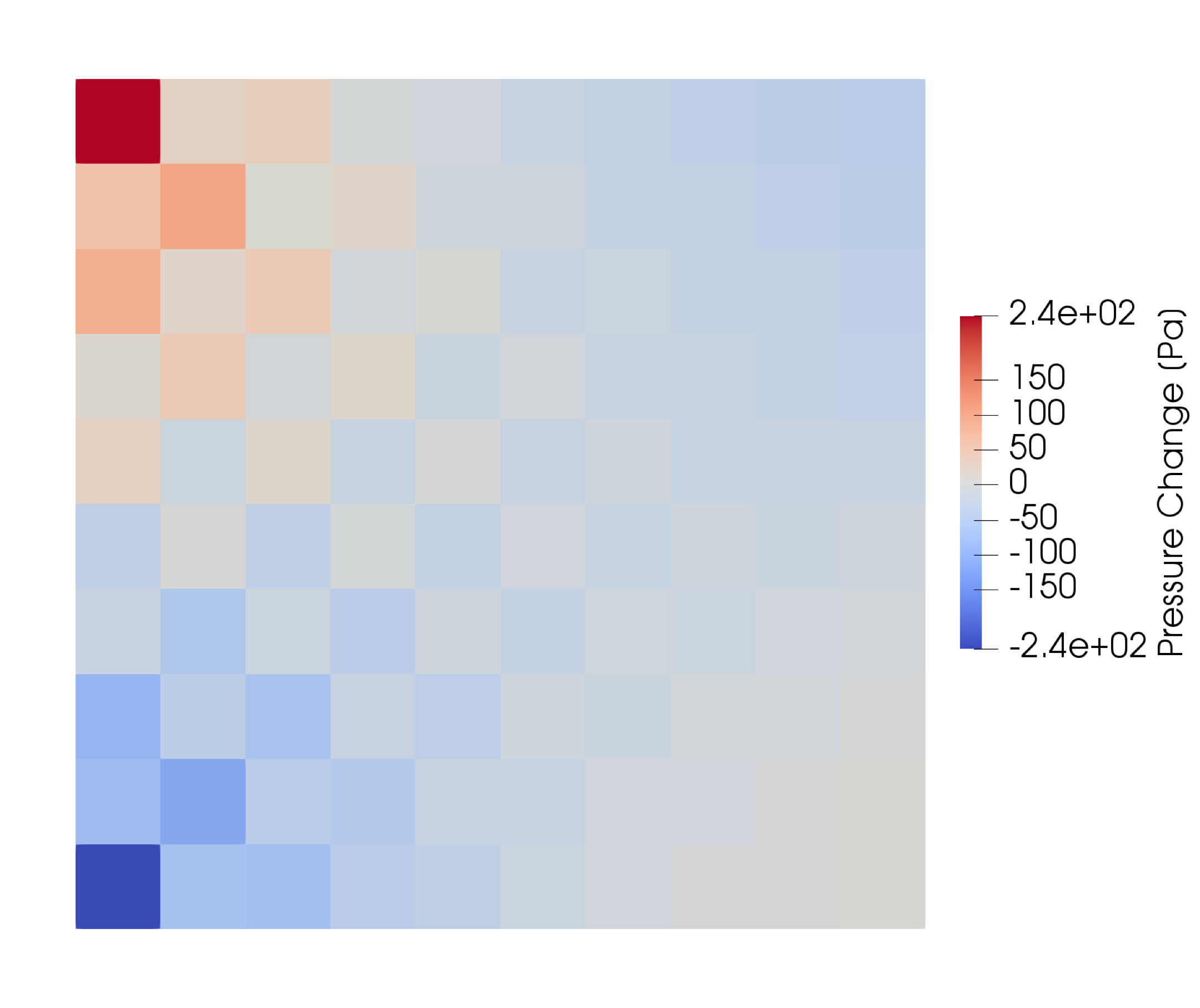}}\\
\caption{Cantilever: Pressure field at various time steps obtained with fully implicit method}
\label{fig:cant_fim}
\end{figure}

\cref{fig:cant_fim} shows the pressure field evolution through one period of loading obtained with the fully implicit method and a time step of 1 day, for reference.
The pressure field is shown after the third time step (which approximately corresponds to the time of maximum loading in the negative $z$ direction) and after the seventh step (which approximately corresponds to the time of maximum loading in the positive $z$ direction).
Clearly the results are polluted by the spurious checkerboard mode, which appears on the first step and persists throughout the simulation.

\begin{figure}
\centering
\subfloat[$t = 3$ days]{\label{sfig:cant_seq_1daydt_3}\includegraphics[width=.4\textwidth]{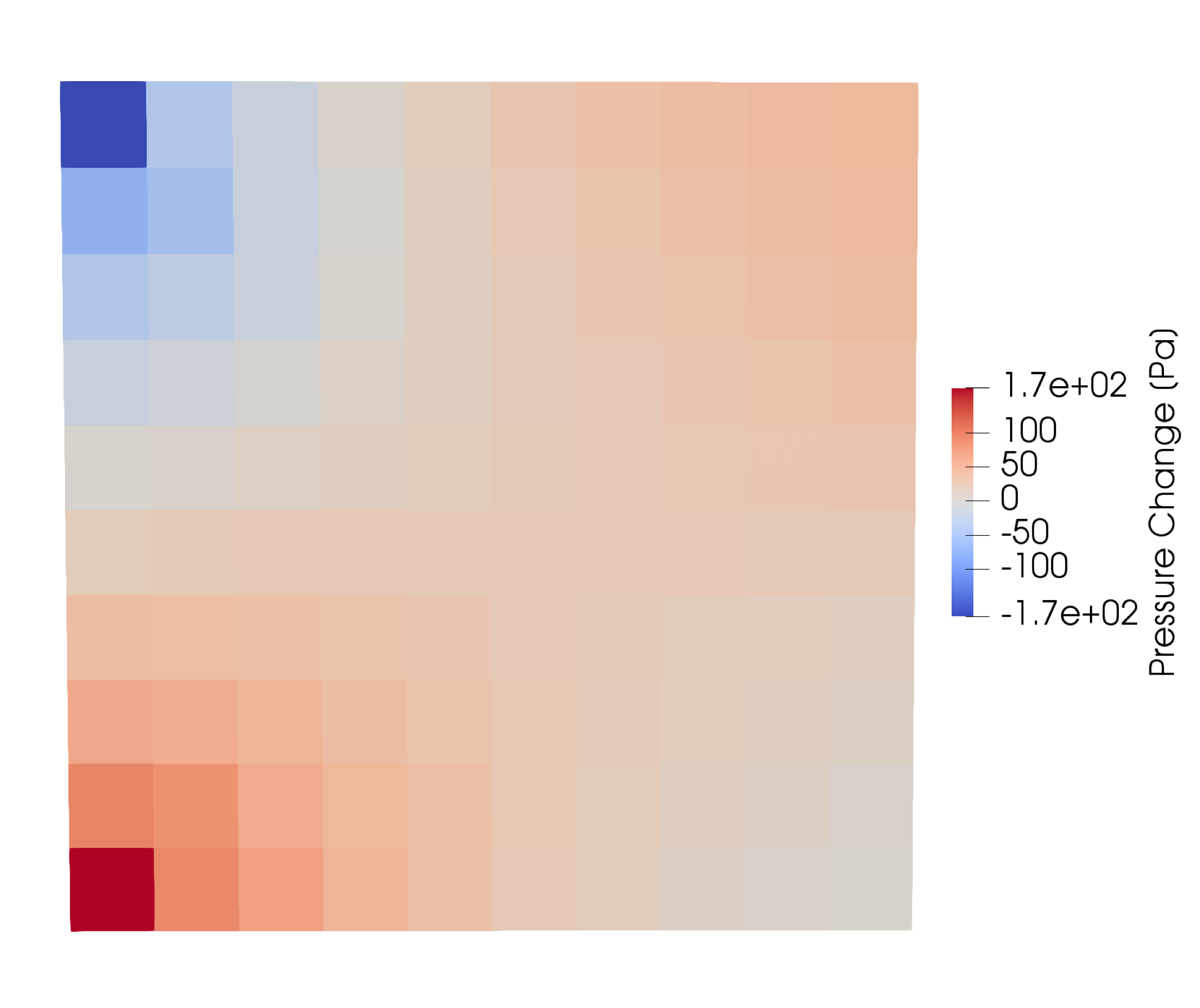}}\hfill
\subfloat[$t = 7$ days]{\label{sfig:cant_seq_1daydt_7}\includegraphics[width=.4\textwidth]{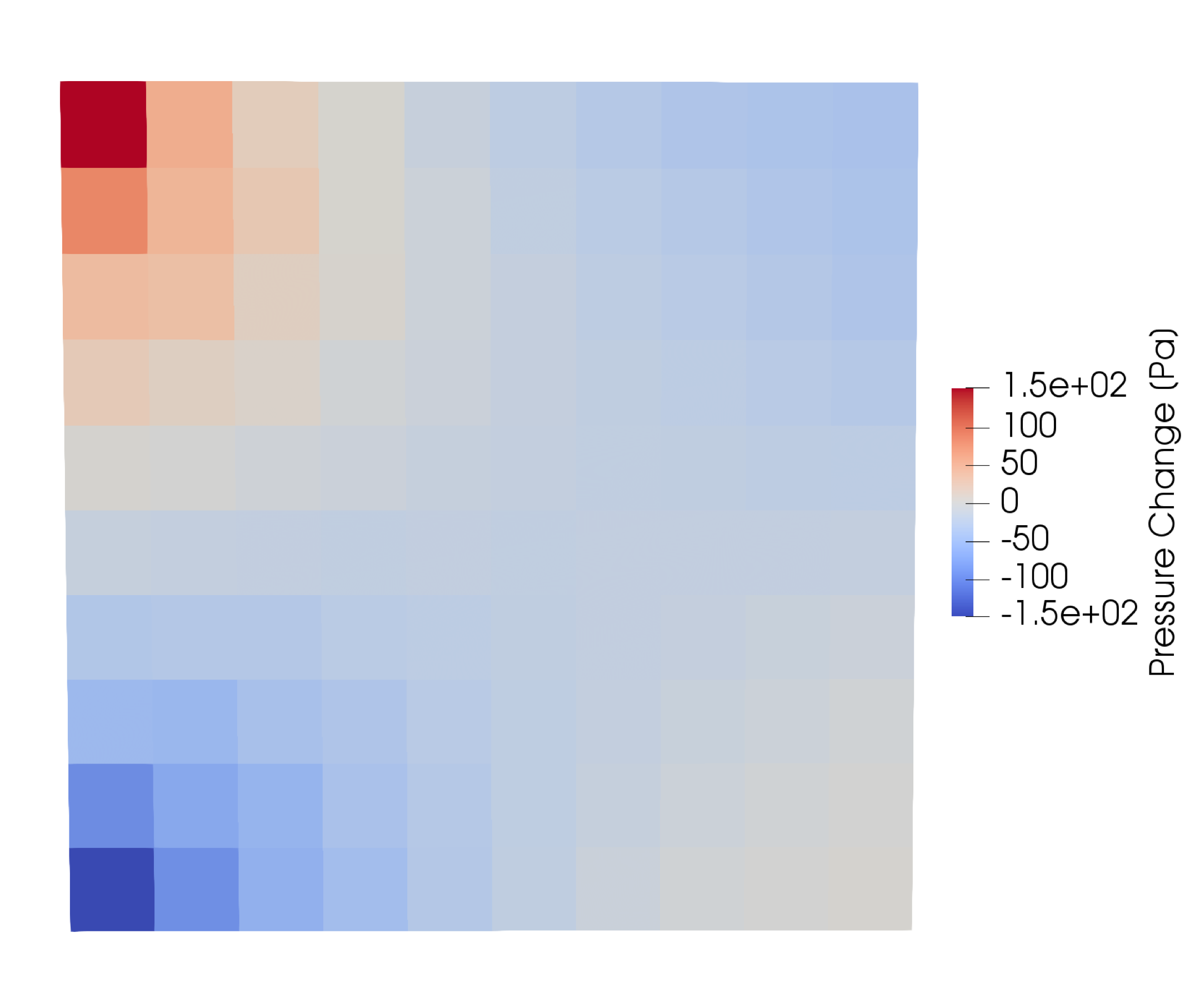}}\\
\caption{Cantilever: Pressure field at various time steps obtained with explicit fixed-stress method and time step of 1 day}
\label{fig:cant_seq_1daydt}
\end{figure}

\begin{figure}
\centering
\subfloat[$t = 3$ days]{\label{sfig:cant_seq_1e-1daydt_3}\includegraphics[width=.4\textwidth]{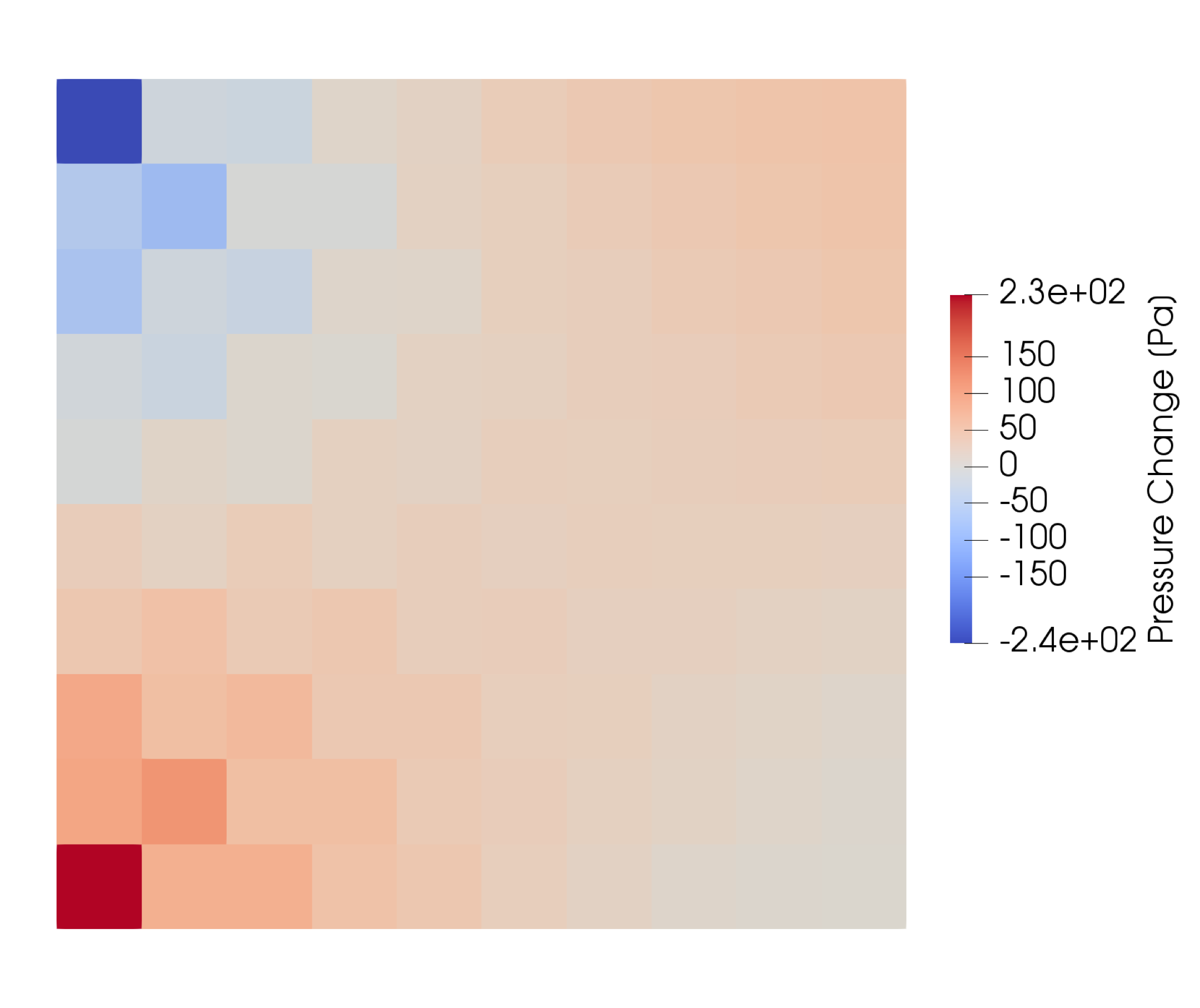}}\hfill
\subfloat[$t = 7$ days]{\label{sfig:cant_seq_1e-1daydt_7}\includegraphics[width=.4\textwidth]{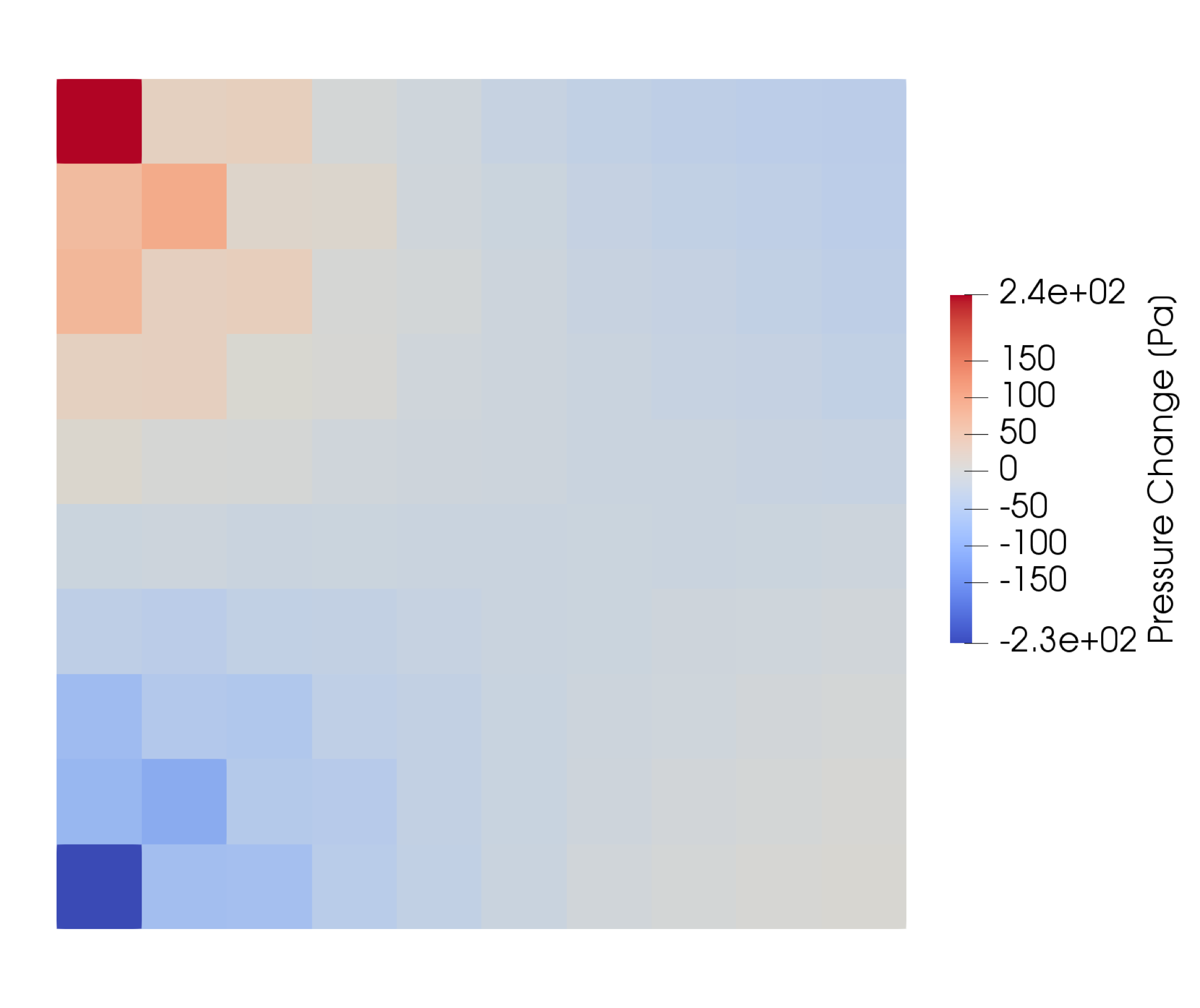}}\\
\caption{Cantilever: Pressure field at various time steps obtained with explicit fixed-stress method and time step of 0.1 day}
\label{fig:cant_seq_1e-1daydt}
\end{figure}

\begin{figure}
\centering
\subfloat[$t = 3$ days]{\label{sfig:cant_seq_1e-2daydt_3}\includegraphics[width=.4\textwidth]{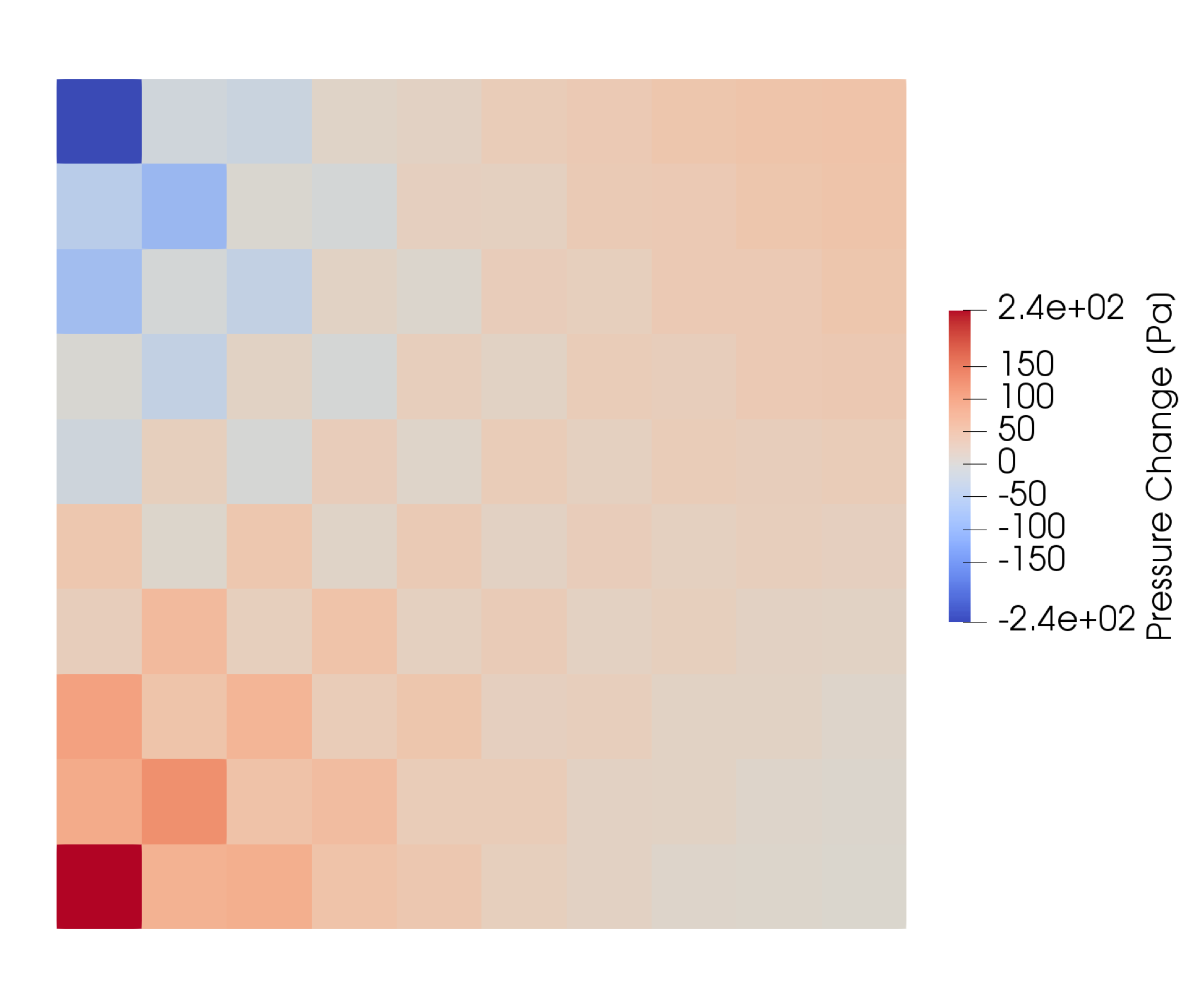}}\hfill
\subfloat[$t = 7$ days]{\label{sfig:cant_seq_1e-2daydt_7}\includegraphics[width=.4\textwidth]{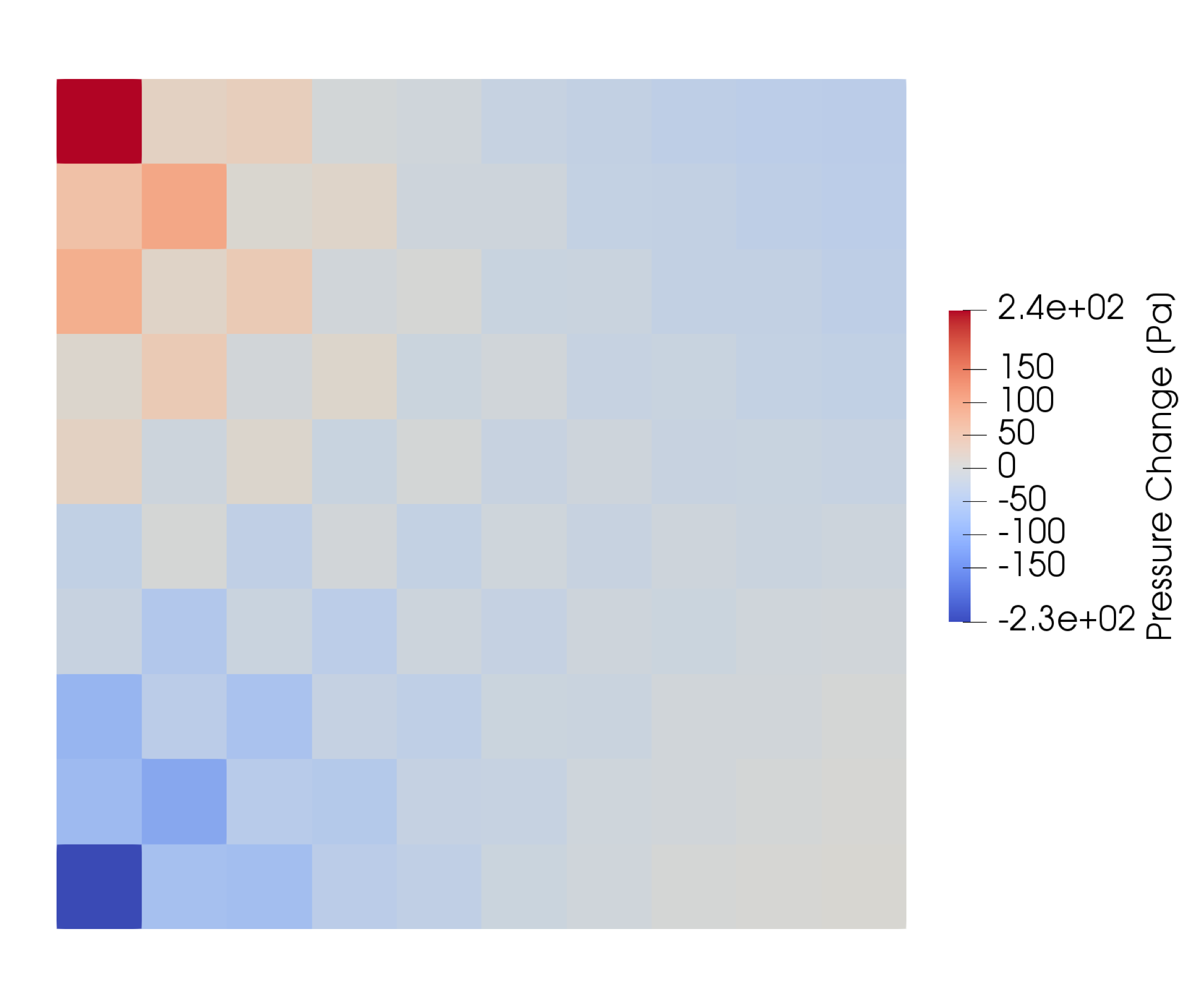}}\\
\caption{Cantilever: Pressure field at various time steps obtained with explicit fixed-stress method and time step of 0.01 day}
\label{fig:cant_seq_1e-2daydt}
\end{figure}

\cref{fig:cant_seq_1daydt} shows the same pressure fields as in \cref{fig:cant_fim}, but obtained with the explicit fixed-stress method with the same 1 day time step.
We see that the pressure field is smooth and does not seem to contain any spurious oscillations.
This confirms our assertion that, for transient problems with large enough time steps, the explicit fixed-stress approach can have a pressure stabilizing effect.
We also wish to study this effect as the time step size decreases.
\cref{fig:cant_seq_1e-1daydt,fig:cant_seq_1e-2daydt} show the corresponding results obtained with time steps of size 0.1 and 0.01 days, respectively.
We see that the checkerboard mode is present in results obtained with a 0.1~day time step, though it is slightly smoothed when compared to the fully implicit results.
When the time step size is further reduced to 0.01 days, the results essentially match the fully implicit case, including the spurious oscillations.
We interpret this as follows: as the time step size is decreased, the solution change during any time step also decreases, and the incompressibility constraints are better enforced and satisfied, which leads to the re-appearance of the spurious pressure modes.

\begin{figure}
\centering
\subfloat[$t = 100$ day]{\label{sfig:cant_seq_1daydt_100}\includegraphics[width=.4\textwidth]{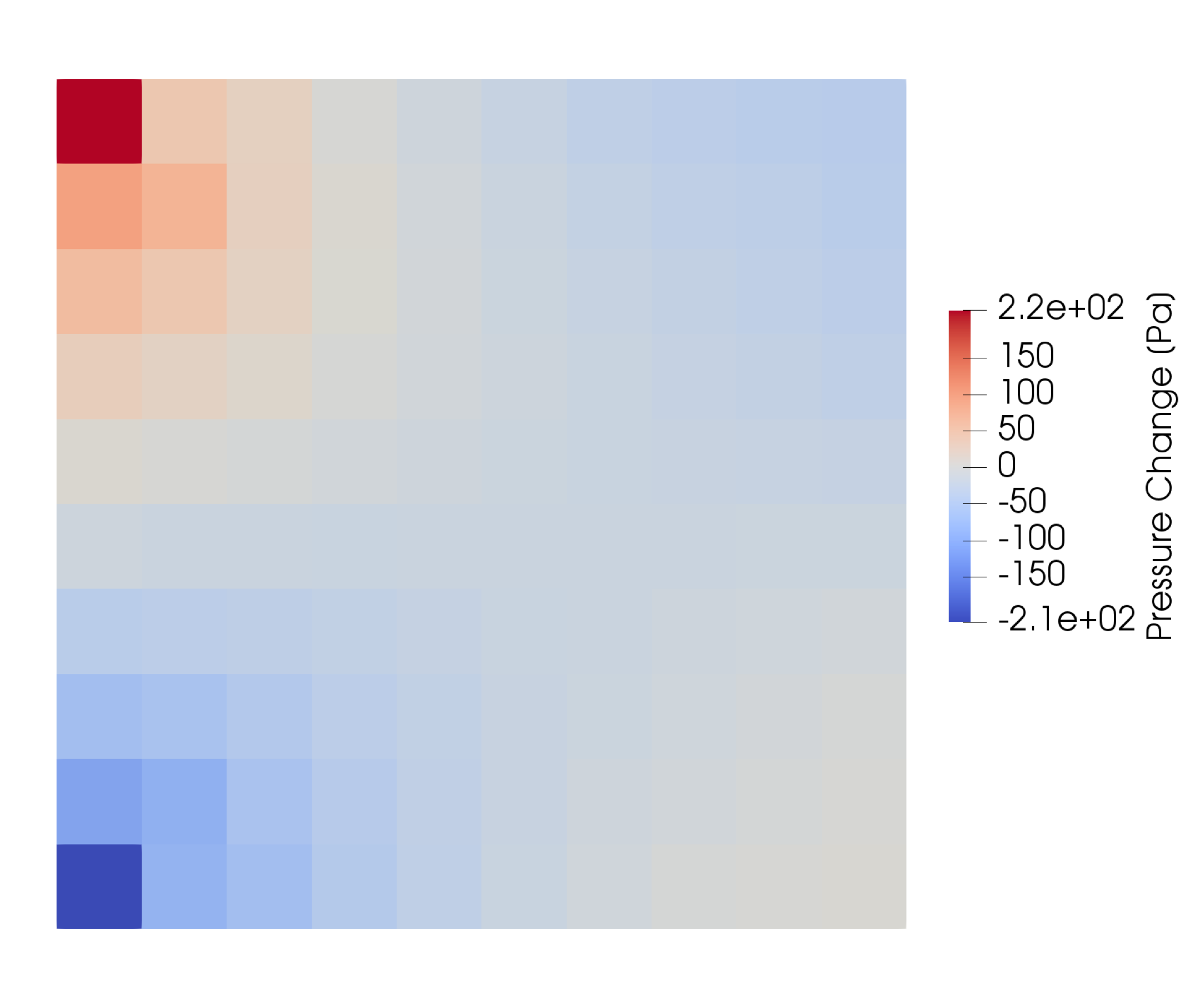}}\hfill
\subfloat[$t = 110$ days]{\label{sfig:cant_seq_1daydt_110}\includegraphics[width=.4\textwidth]{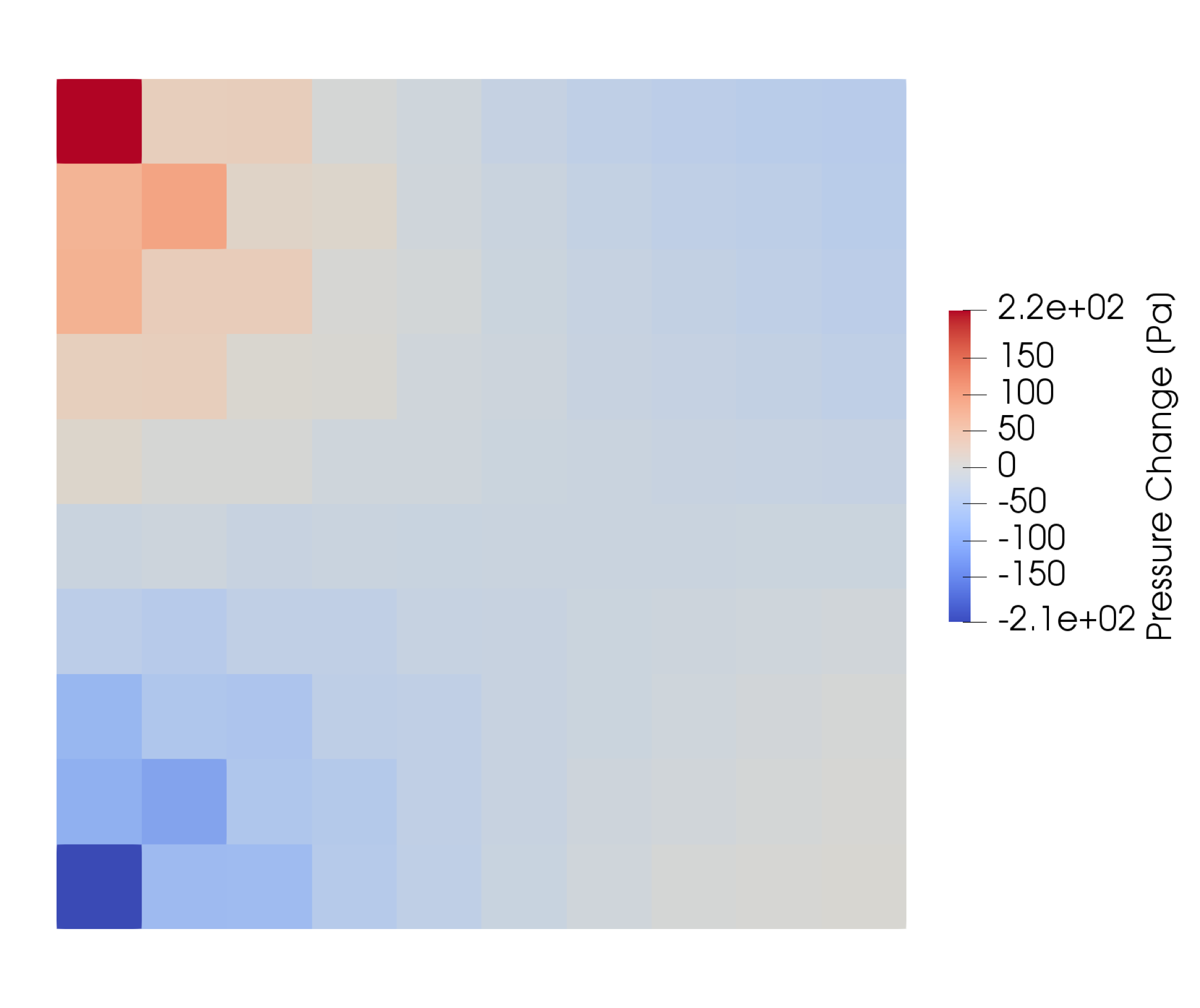}}\\
\subfloat[$t = 120$ days]{\label{sfig:cant_seq_1daydt_120}\includegraphics[width=.4\textwidth]{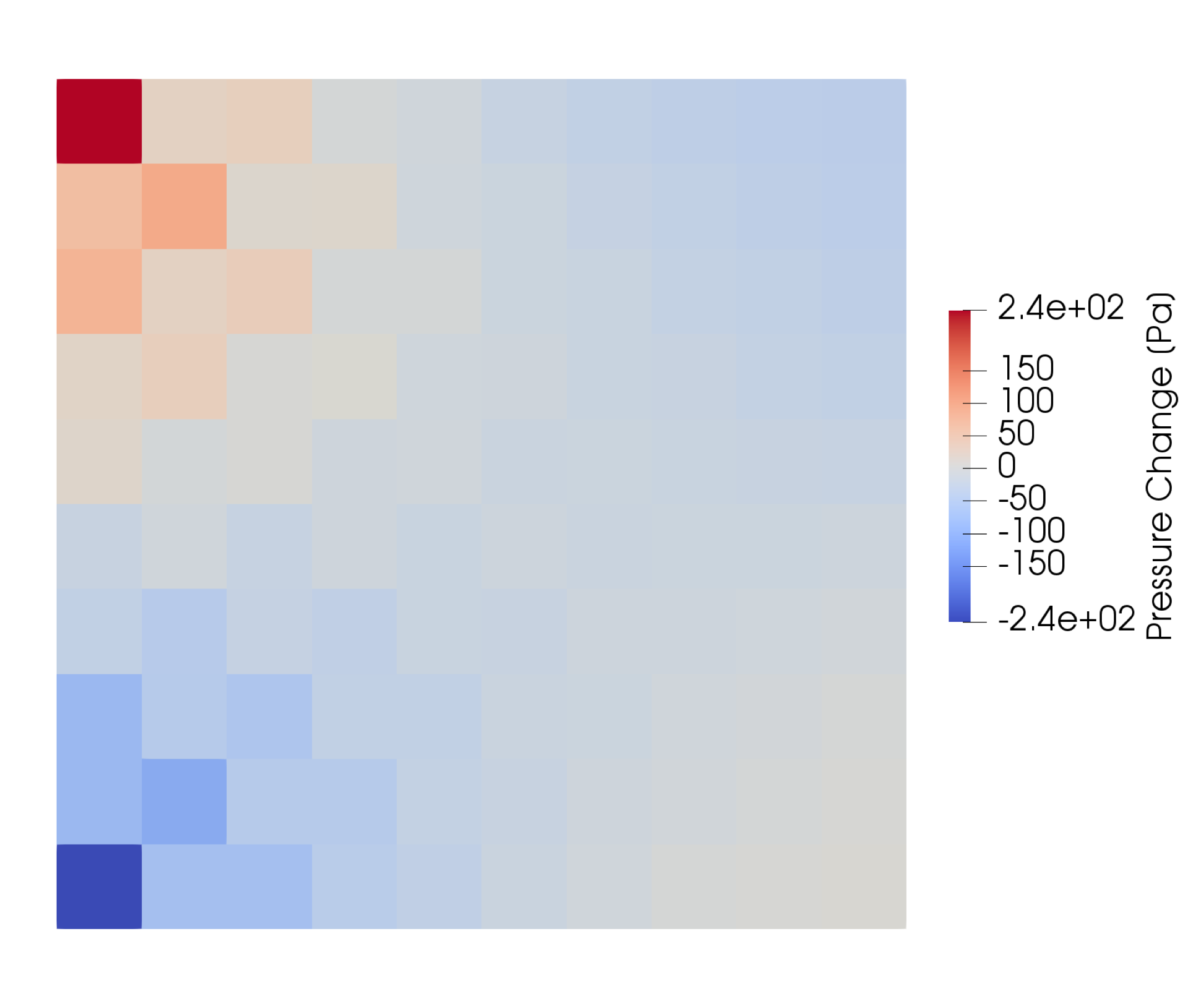}}\hfill
\subfloat[$t = 130$ days]{\label{sfig:cant_seq_1daydt_130}\includegraphics[width=.4\textwidth]{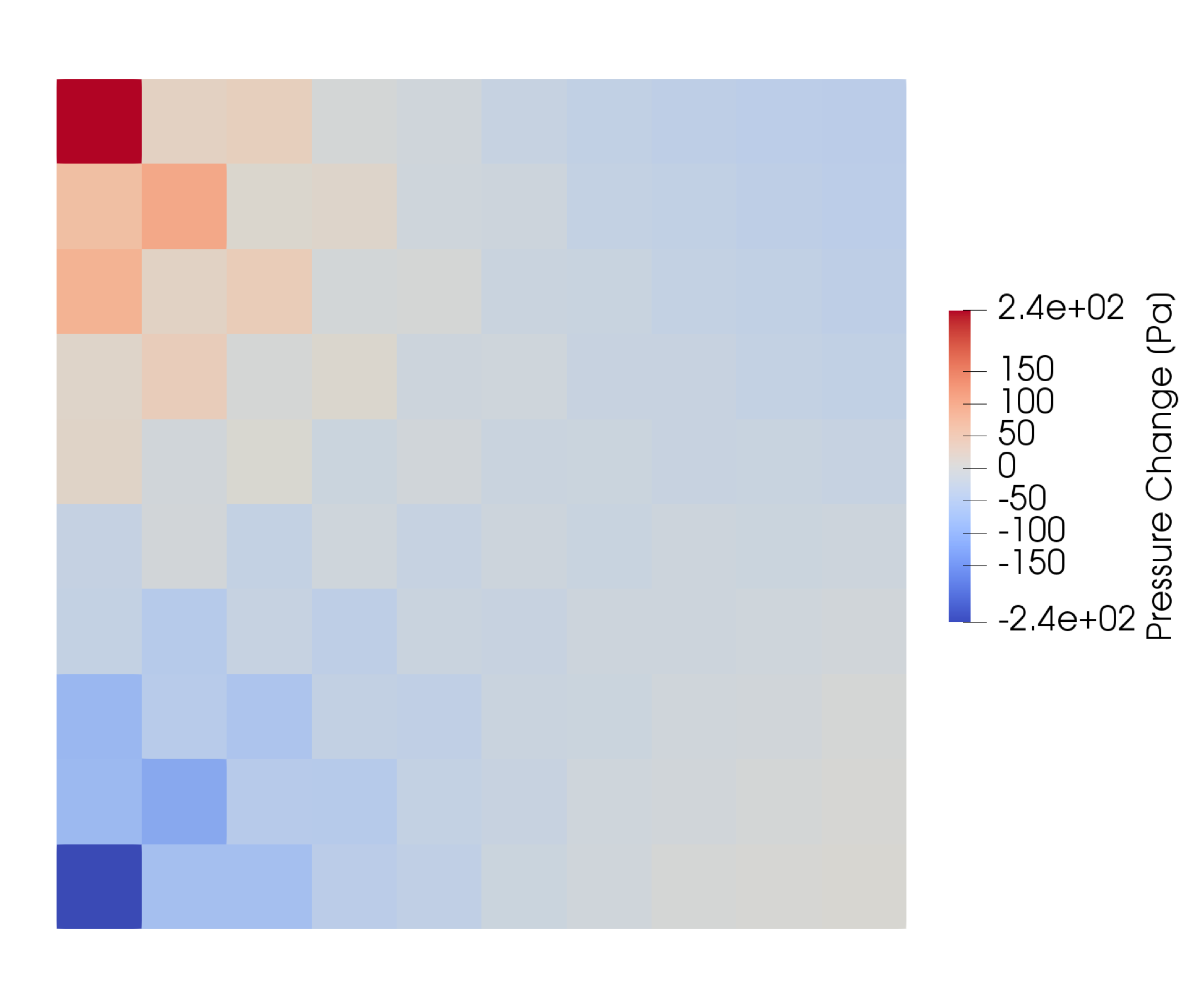}}\\
\caption{Cantilever: Pressure field at various time steps obtained with explicit fixed-stress method and time step of 1 day with forcing frozen at maximum value}
\label{fig:cant_seq_1daydt_steady}
\end{figure}

As a last numerical test, we consider the behavior of the explicit scheme as steady state is approached when using the 1 day time step which smoothed the spurious pressure modes in the transient simulations above.
We continue applying the sinusoidal load for 100 days with the explicit fixed-stress method, and the resulting pressure field is given by \cref{sfig:cant_seq_1daydt_100}.
Clearly the spurious pressure mode is still not visible at this point.
At $t = 100$ days, we freeze the forcing at its maximum value in the positive $z$ direction and continue time stepping.
The problem is now steady, and we plot the pressure field obtained after 10, 20, and 30 extra time steps in the remainder of \cref{fig:cant_seq_1daydt_steady}.
Even though the spurious oscillations remained invisible during the transient phase of the simulation, we see that continuing to advance in time at steady state results in a~solution which converges to the oscillatory solution obtained by the fully implicit method. 

\subsection{Staircase}

\begin{figure}
\centering
\includegraphics[width=0.4\textwidth]{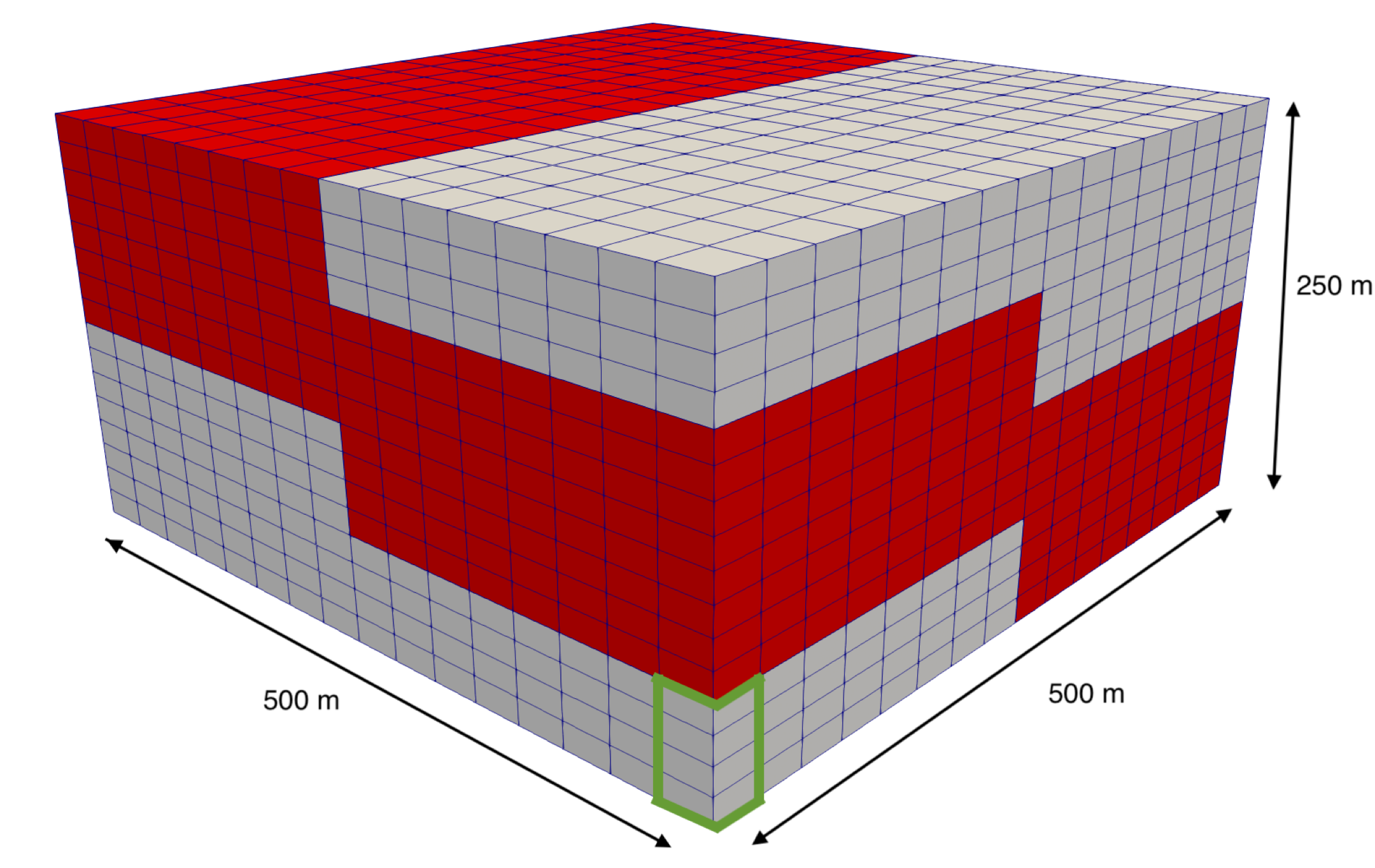}
\caption{Staircase problem setup}
\label{fig:stair_setup}
\end{figure}

Our second numerical example consists of a high-permeability channel region and a low permeability burden region coupled together in a spiral staircase pattern, as seen in \cref{fig:stair_setup}.
This problem has been studied previously in \cite{white2008stabilized, camargo2021macroelement, aronson2024FS}, and is a useful example as it remains simple but is perhaps more closely related to what one would expect in CO$_2$ sequestration, namely the undrained conditions occur only in the burden region while injection occurs into the high-permeability channel region.
We again consider single-phase poromechanics for simplicity, and injection is performed into the cells marked in green in \cref{fig:stair_setup} at a rate of 1 kg/s for 30~years.
The skeleton in both the channel and barrier regions is modeled as isotropic, linear elastic with bulk modulus of 5 GPa and Poisson ration of 0.25.
Roller boundary conditions are applied on all sides except for the top, which is left free.

Within the channel region (shown in gray in \cref{fig:stair_setup}), the skeleton permeability is set to $9.8\cdot 10^{-13}$~m$^2$ and the porosity is set to 0.2, while in the barrier region (shown in red in \cref{fig:stair_setup}), the permeability is set to zero and the porosity is set to 0.05.
Thus the problem is drained in the channel and undrained in the barrier.
The fluid is again assumed to be fully incompressible with density 1000 kg/m$^3$, and no flow boundary conditions are imposed on all sides of the domain.

\begin{figure}
\centering
\subfloat[$t = 15$ years]{\label{sfig:stair_fim_15}\includegraphics[width=.4\textwidth]{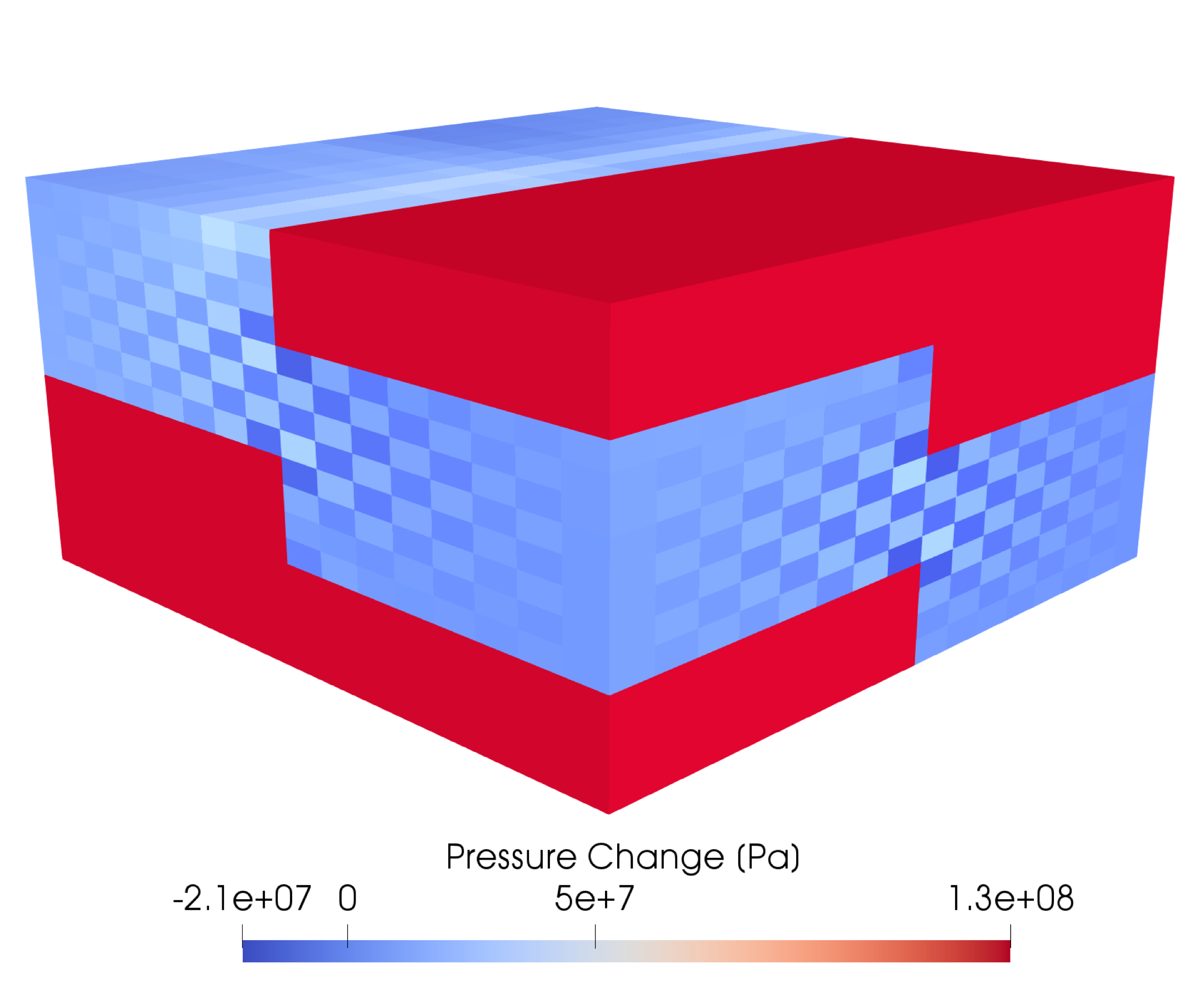}}\hfill
\subfloat[$t = 30$ years]{\label{sfig:stair_fim_30}\includegraphics[width=.4\textwidth]{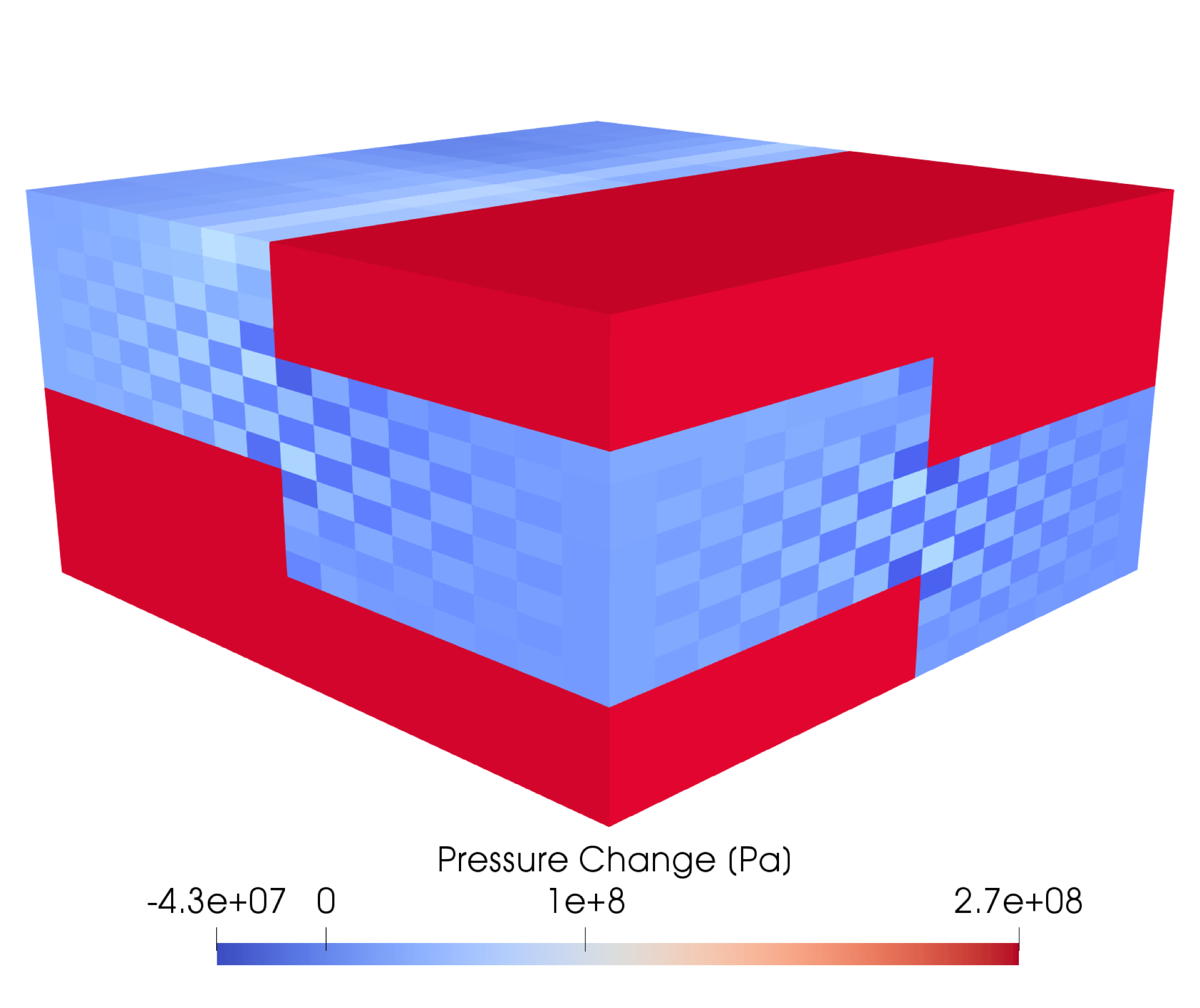}}\\
\caption{Staircase: Pressure field at various time steps obtained with fully implicit method}
\label{fig:stair_fim}
\end{figure}

Like the previous example, we first show a reference result obtained with a fully implicit method.
\cref{fig:stair_fim} shows the pressure field after 15 and 30 years with a time step of 1 month.
There are strong pressure oscillations in the undrained burden region, and we note that, though it is not shown here, they persist when the time step size is increased.

\begin{figure}
\centering
\subfloat[$t = 15$ years]{\label{sfig:stair_seq1month_15}\includegraphics[width=.4\textwidth]{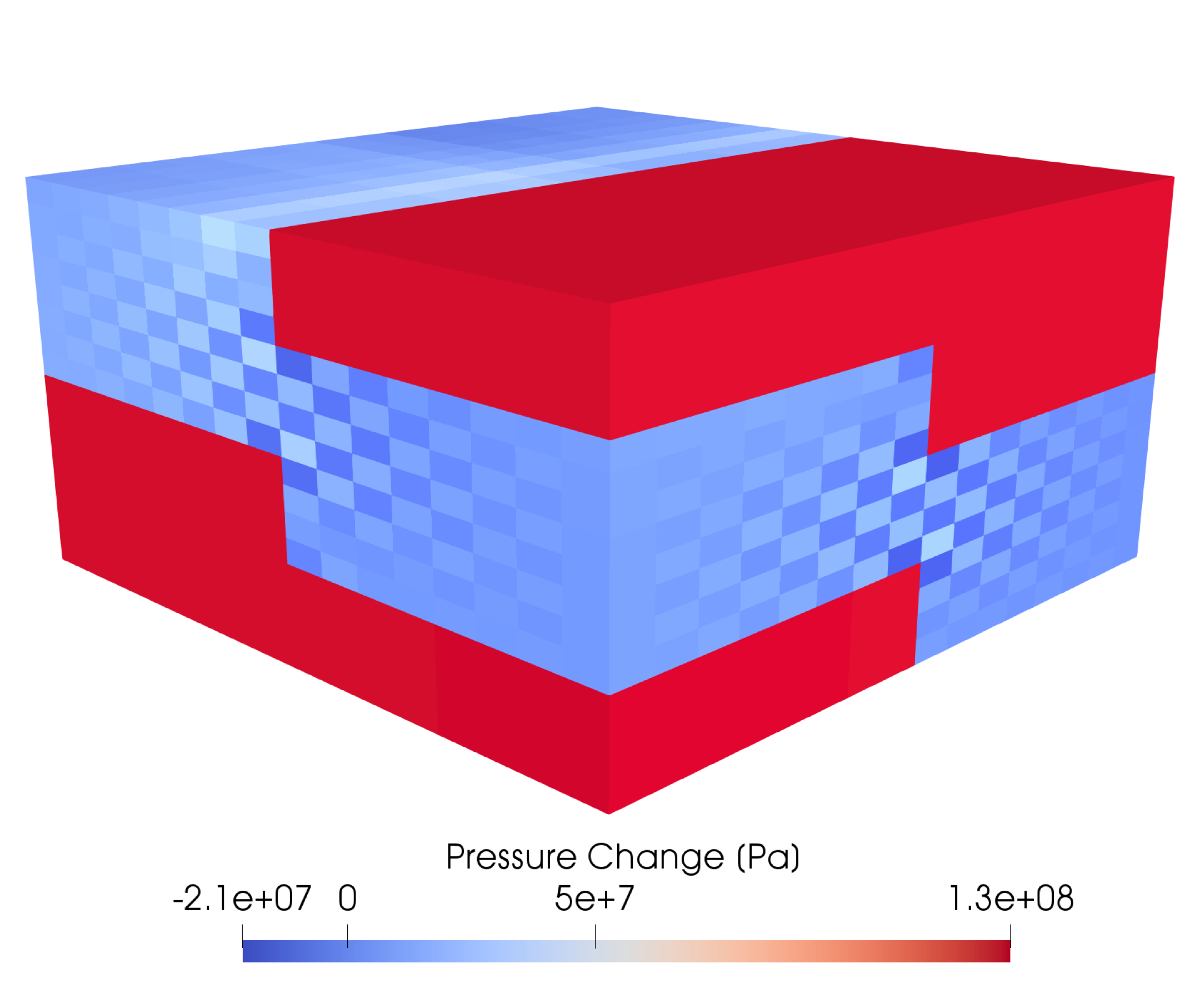}}\hfill
\subfloat[$t = 30$ years]{\label{sfig:stair_seq1month_30}\includegraphics[width=.4\textwidth]{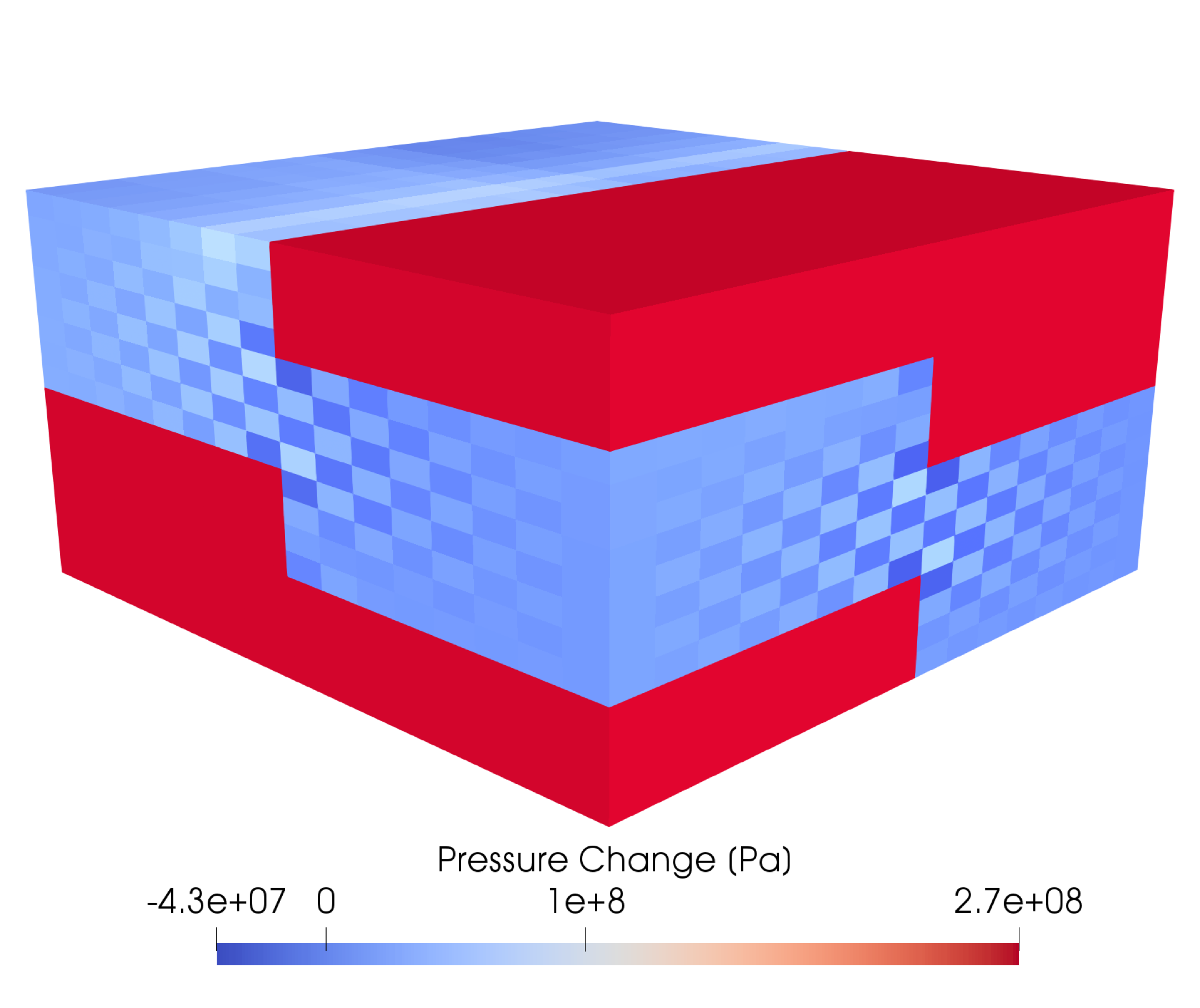}}\\
\caption{Staircase: Pressure field at various time steps obtained with explicit fixed-stress method and time step of 1 month}
\label{fig:stair_seq_1month}
\end{figure}

\begin{figure}
\centering
\subfloat[$t = 15$ years]{\label{sfig:stair_seq1year_15}\includegraphics[width=.4\textwidth]{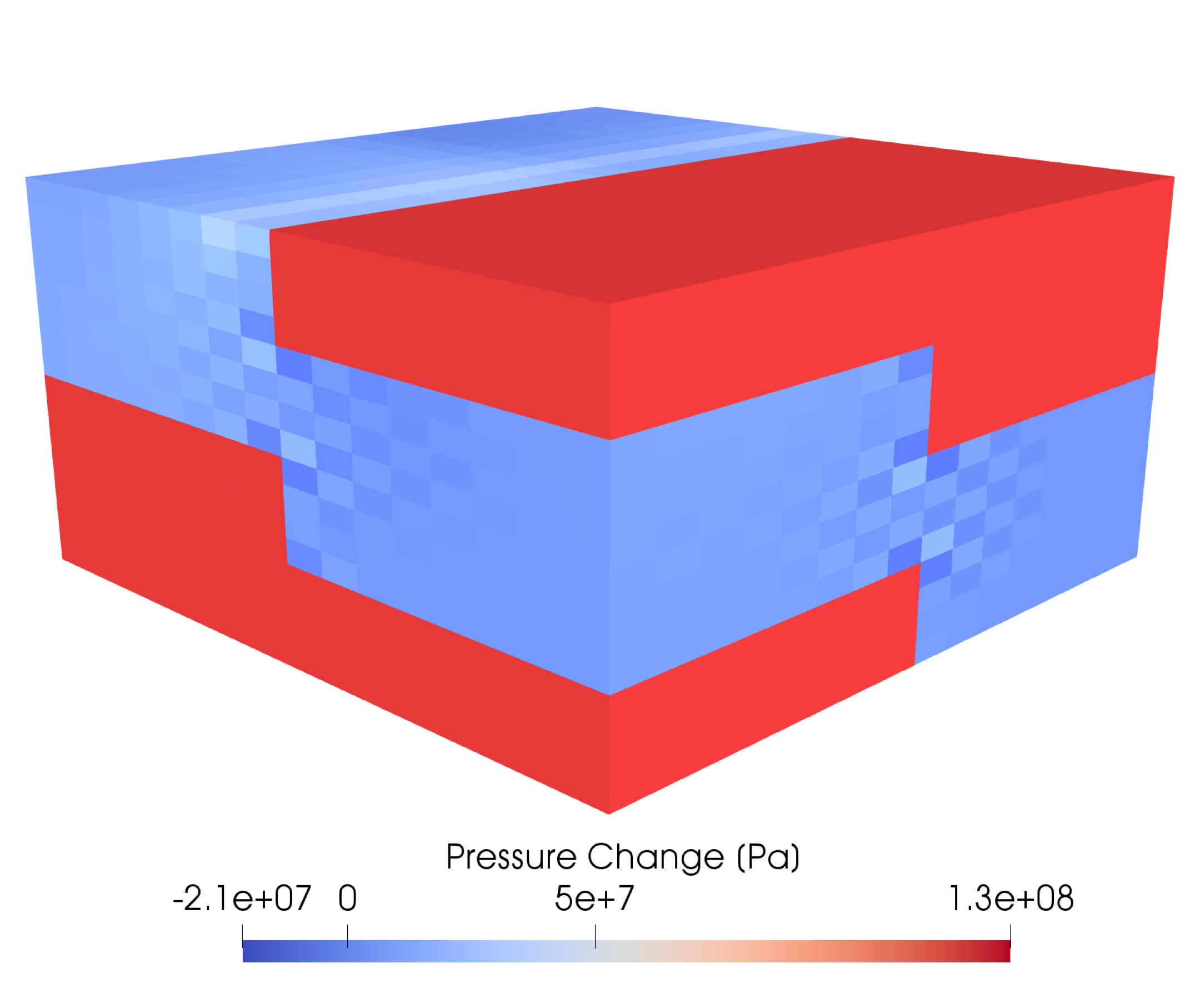}}\hfill
\subfloat[$t = 30$ years]{\label{sfig:stair_seq1year_30}\includegraphics[width=.4\textwidth]{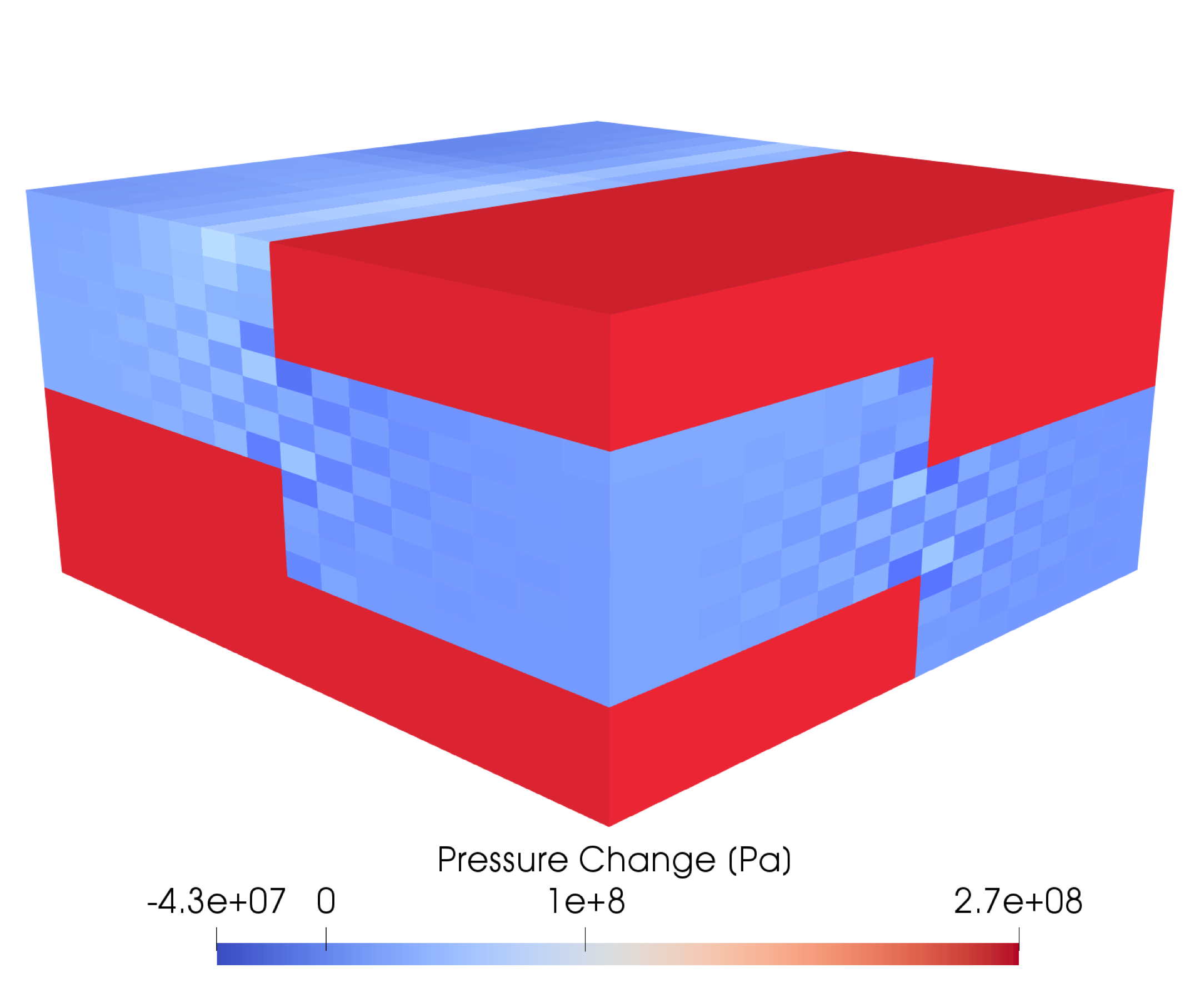}}\\
\caption{Staircase: Pressure field at various time steps obtained with explicit fixed-stress method and time step of 1 year}
\label{fig:stair_seq_1year}
\end{figure}

\begin{figure}
\centering
\subfloat[$t = 15$ years]{\label{sfig:stair_seq5year_15}\includegraphics[width=.4\textwidth]{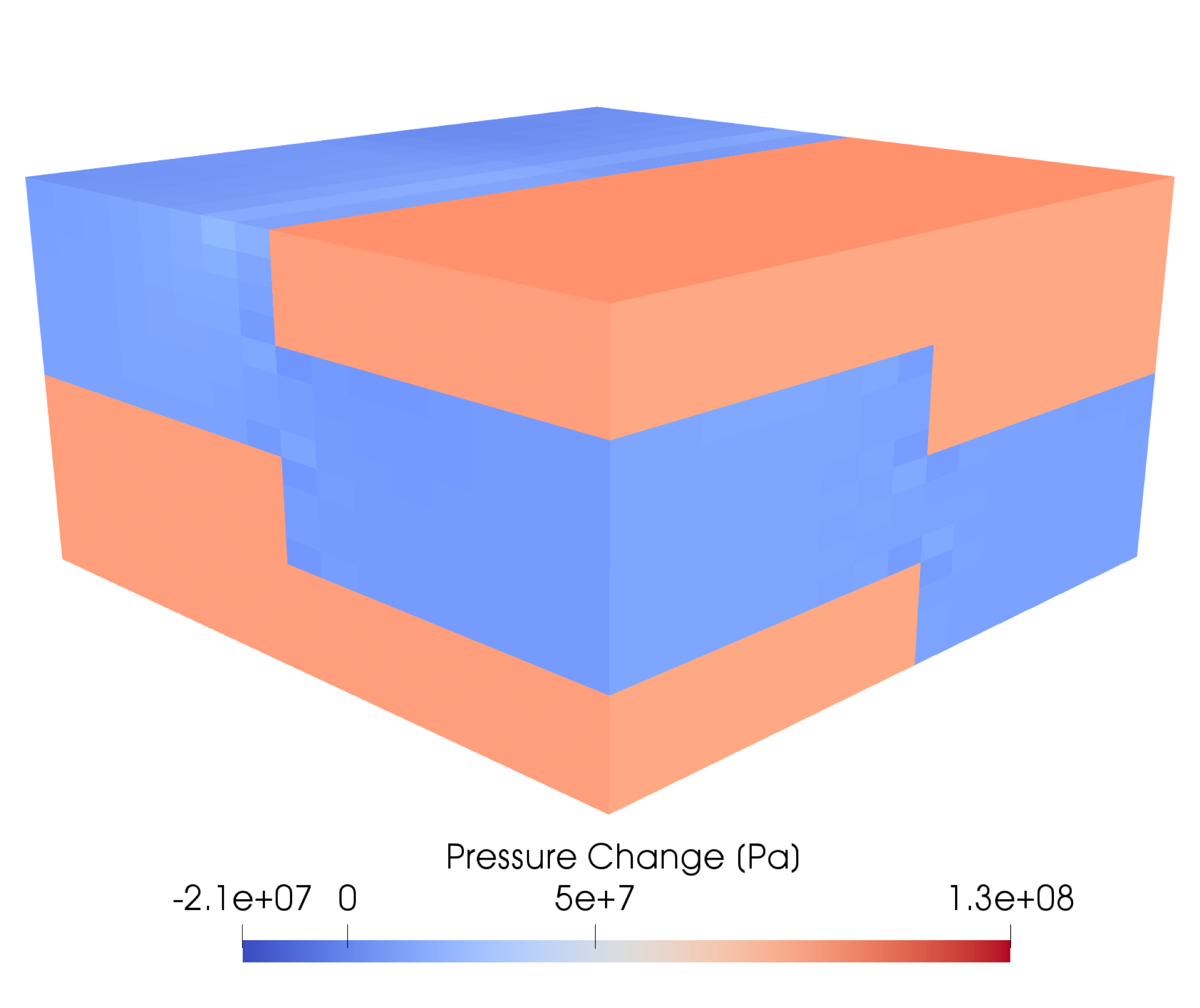}}\hfill
\subfloat[$t = 30$ years]{\label{sfig:stair_seq5year_30}\includegraphics[width=.4\textwidth]{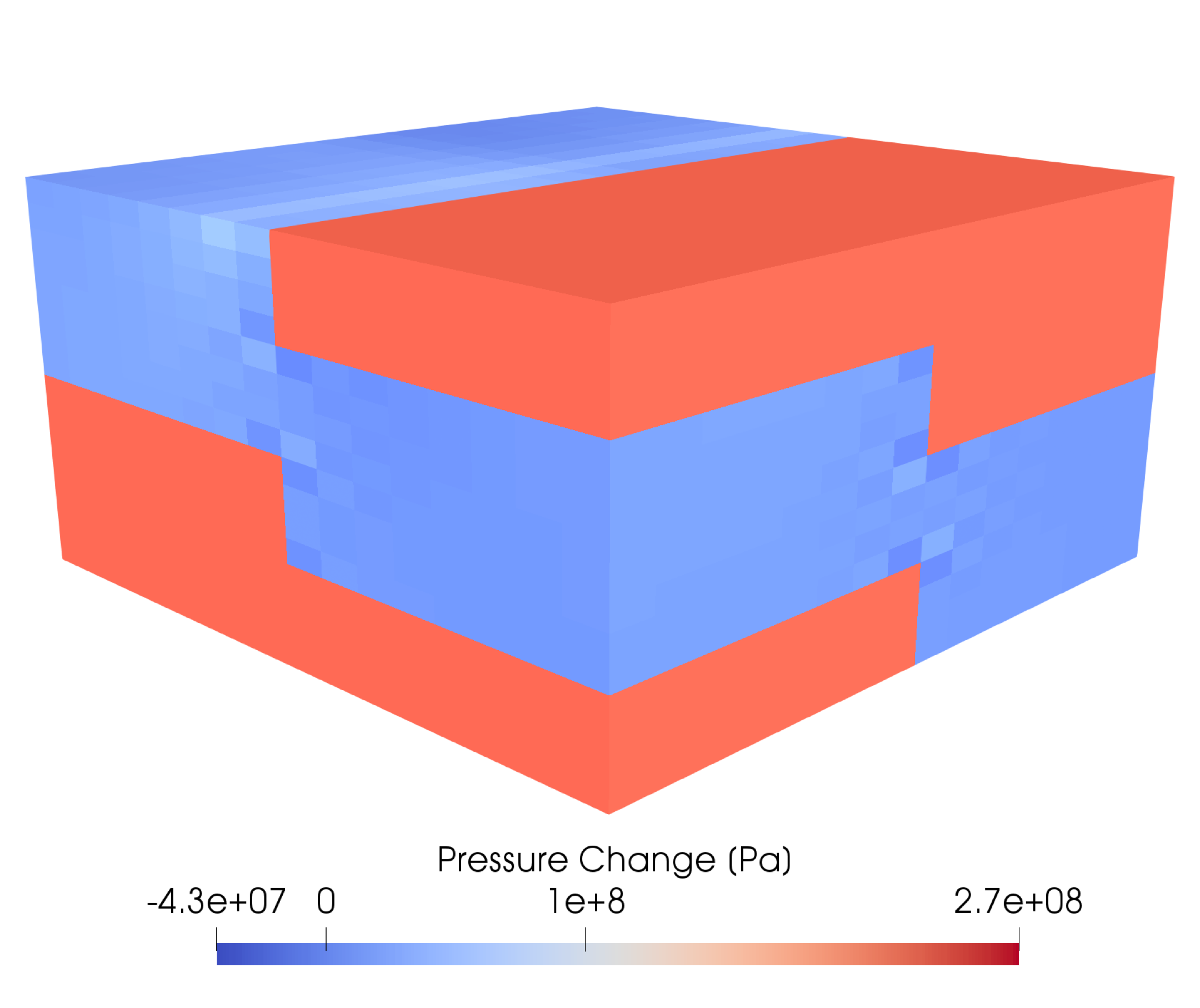}}\\
\caption{Staircase: Pressure field at various time steps obtained with explicit fixed-stress method and time step of 5 years}
\label{fig:stair_seq_5year}
\end{figure}

\cref{fig:stair_seq_1month,fig:stair_seq_1year,fig:stair_seq_5year} show the same results using the explicit fixed-stress method and time steps of 1~month, 1~year, and 5~years, respectively.
We see the same trends as in the cantilever example -- when the time step size is set to be 1 month, the solution is almost identical to the fully implicit solution, and there is essentially no pressure stabilizing effect.
When the time step is increased to 1 year, we do see a small about of pressure smoothing, and further increasing to 5 years removes most of the spurious oscillations.
However, we also see that increasing the time step size in the explicit scheme also introduced noticeable error into the solution within the channel region.
This is expected, but highlights the dilemma that appears when attempting to use an explicit coupling scheme to stabilize the pressure in this kind of example.
In~particular, since the pressure solution in the low permeability barrier region changes very slowly, a large time step is needed for the transient effects to damp oscillatory modes.
The channel solution, on the other hand, experiences much larger changes and thus is subjected to the larger transient errors that arise from a large time step. In this region, however, these effects are not necessary to counteract any spurious modes, and so the damping is unwanted.
Finally, while we do not show it here for brevity, we remark that at steady state the behavior of this example is the same as the cantilever example.
In particular, if we were to continue time stepping beyond 30 years with no injection, the solutions obtained with the explicit fixed-stress scheme with any time step size will converge to a solution with the same pressure oscillations in the barrier as obtained with a fully implicit method.

\FloatBarrier

\subsection{High Island 24L}

Our final numerical example is a realistic setup of field-scale CO$_2$ sequestration.
We consider injection into a geologic model of High Island 24L field (HI24L), located in the Gulf of Mexico.
This region has been extensively studied and characterized \cite{ruiz2019characterization}, and is an attractive candidate for CO$_2$ storage.

\begin{figure}
\centering
\subfloat[Full simulation mesh]{\label{sfig:HI_mesh}\includegraphics[width=.4\textwidth]{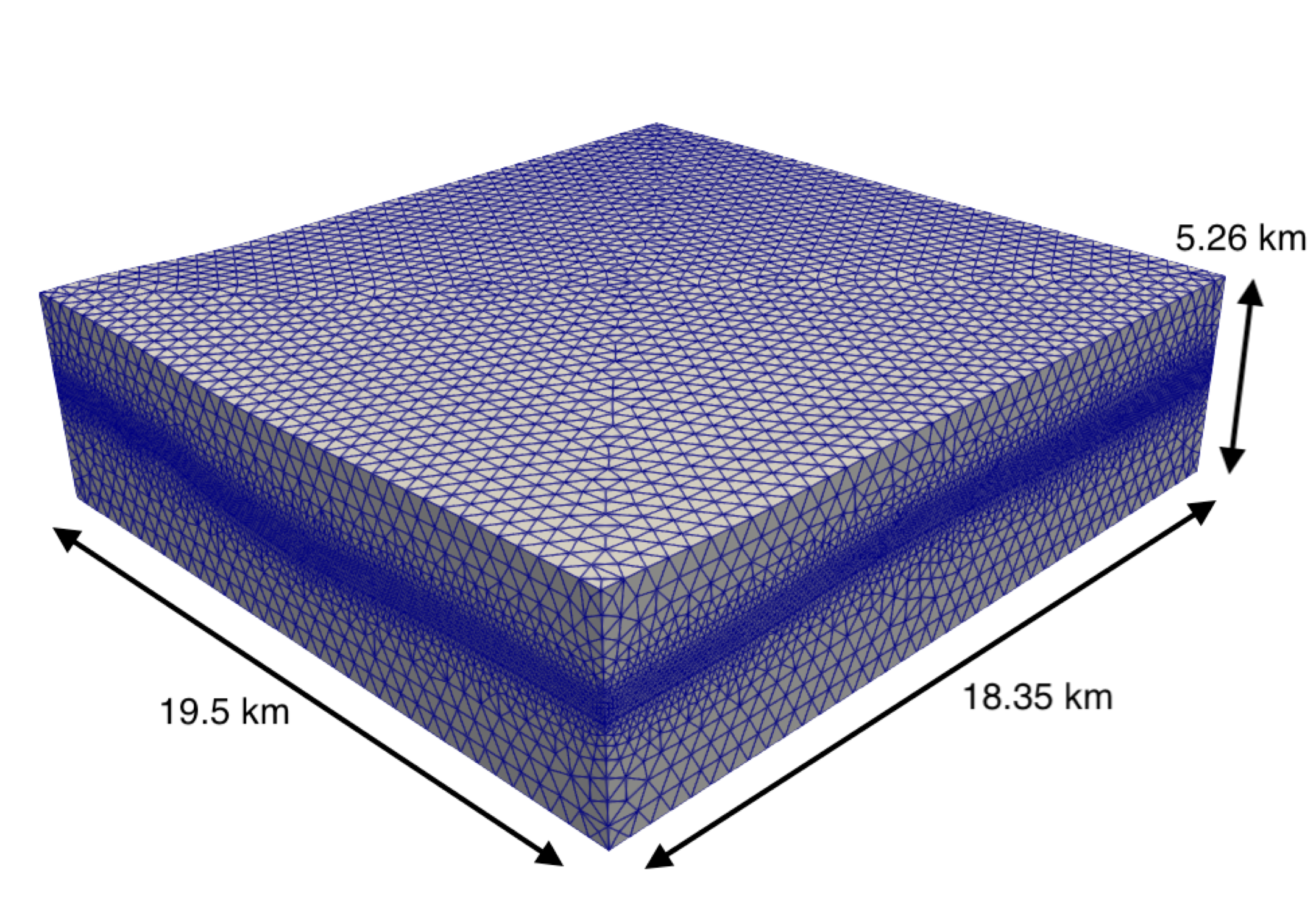}}\hfill
\subfloat[Variability of aquifer permeability and porosity. Permeability varies between  $10^{-16} - 5\cdot 10^{-13}$~m$^2$ while porosity varies between 5 - 25\%.]{\label{sfig:HI_perm}\includegraphics[width=.3\textwidth]{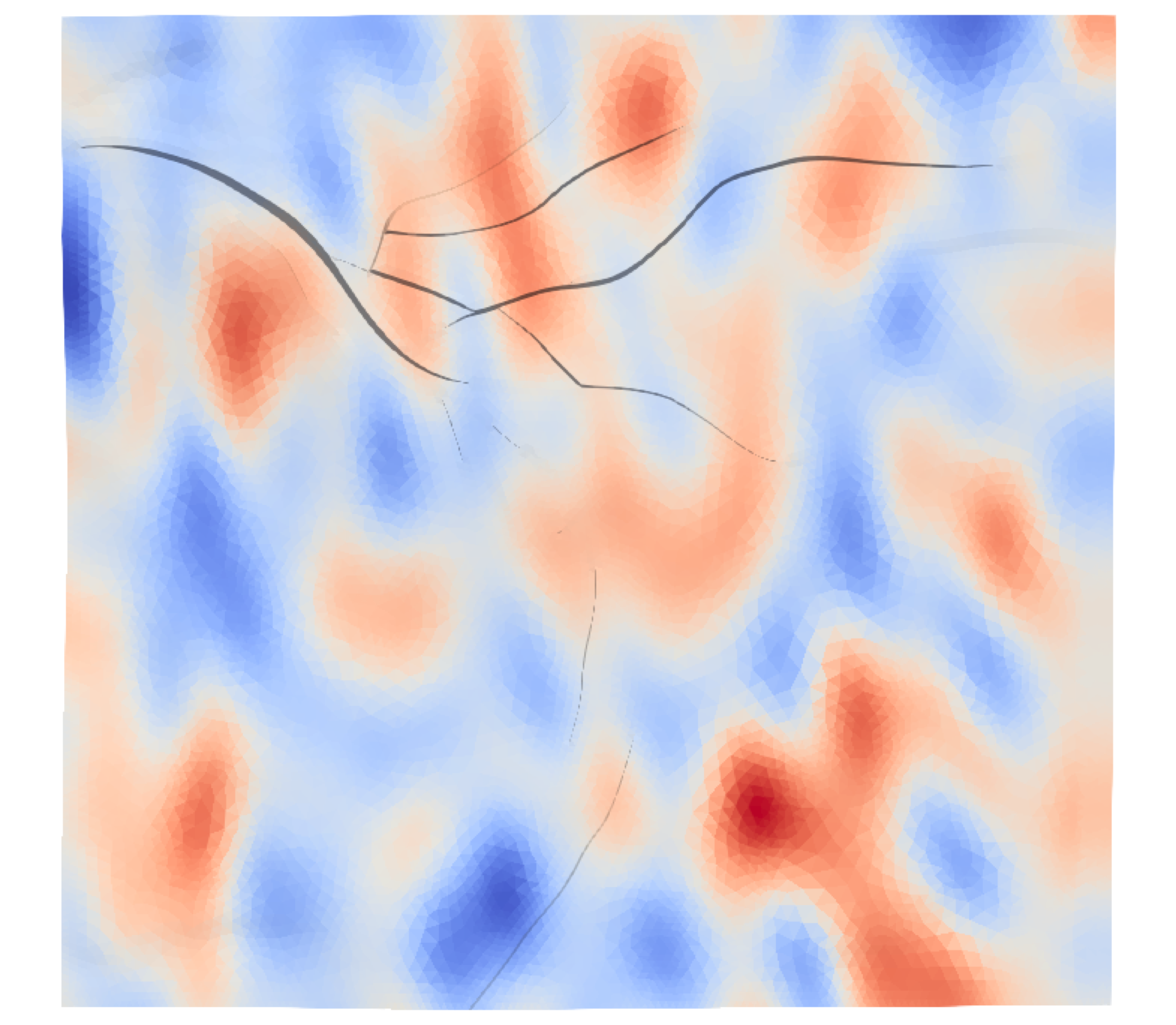}}\\
\caption{HI24L: Problem setup}
\label{fig:HI_setup}
\end{figure}

The simulation mesh is shown in \cref{sfig:HI_mesh}, and consists of a high permeability aquifer sandwiched between low permeability seal units.
Within the aquifer region, the skeleton density is set to 2700 kg/m$^3$, the drained bulk modulus is set to 9.4 GPa, and the Poisson ratio is set to 0.25.
The permeability and porosity fields are heterogeneous within the aquifer, with the permeability being isotropic and varying between $10^{-16}$~m$^2$ and $5\cdot 10^{-13}$~m$^2$ and the porosity varying between 5\% and 25\%.
The profile for both fields is shown in \cref{sfig:HI_perm}.
Within the seal units, the skeleton density remains 2700 kg/m$^3$, the drained bulk modulus is 11.5 GPa, and the Poisson ratio is 0.3.
The permeability is isotropic and homogeneous, with value $10^{-20}$~m$^2$ in each direction.
The porosity is also homogeneous with value 5\%.
The top surface of the skeleton is free, while all other boundaries are defined with roller conditions.
All boundaries prevent flow in or out of the domain.

On the fluid side, the CO$_2$ solubility is defined as in \cite{duan2003improved}, and the CO$_2$-rich phase density is defined using~\cite{span1996new} while the viscosity is determined as in \cite{fenghour1998viscosity}.
The brine phase density and viscosity are determined using \cite{phillips1981technical}.
CO$_2$ is injected at a rate of 90 kg/s in the center of the reservoir for 30 years.

\begin{figure}
\centering
\subfloat[Fully Implicit]{\label{sfig:HI_FIM}\includegraphics[width=.4\textwidth]{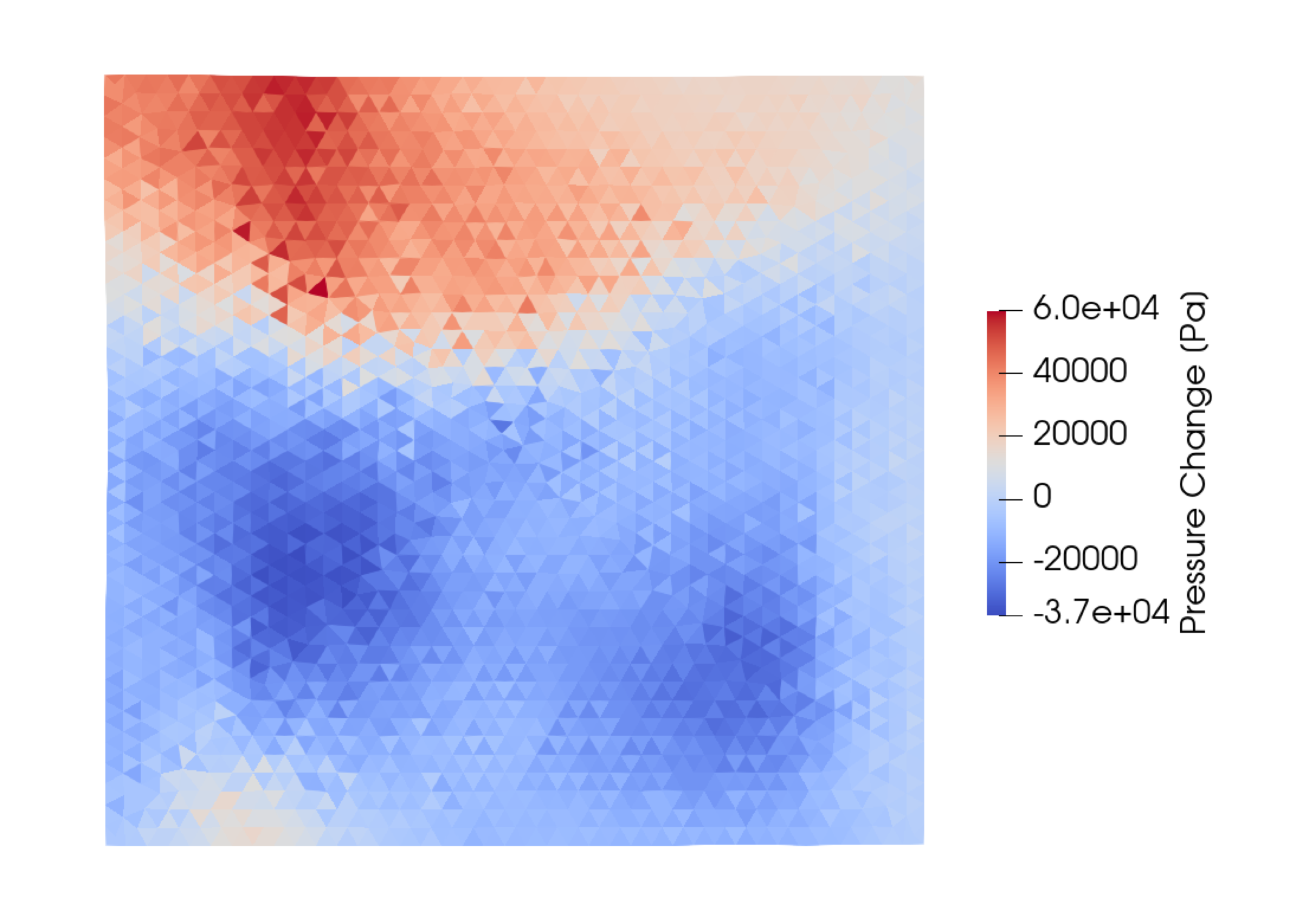}}\hfill
\subfloat[Explicit Fixed-Stress]{\label{sfig:HI_SEQ}\includegraphics[width=.4\textwidth]{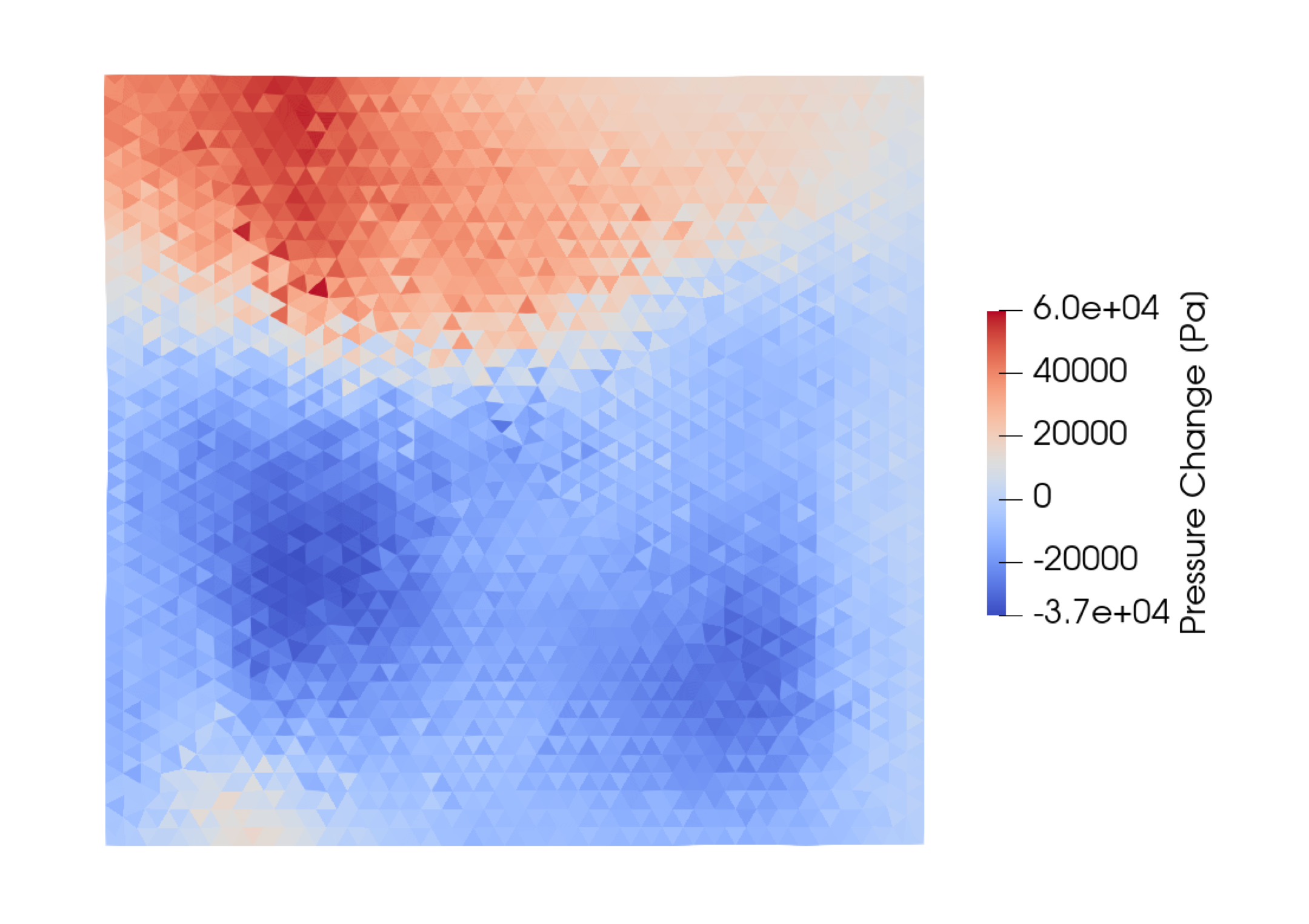}}\\
\caption{HI24L: Seafloor pressure field after 30 years of injection with time step of 6 months}
\label{fig:HI}
\end{figure}

\cref{fig:HI} shows the change in the pressure field at the top of the upper seal unit, which represents the seafloor, using both the fully implicit method and the explicit fixed-stress split, and a time step of size 6~months.
Clearly both are polluted with the same spurious pressure oscillations.
Given the previous results in the staircase section, this is perhaps not surprising and achieving a pressure stabilizing effect from the explicit coupling would likely require a much larger time step.
However, in this case the time step size is now limited by the performance of the nonlinear solution strategy for the flow problem.
For this example, even increasing the time step from 6 months to 1 year prevents the Newton solver from converging in 100 iterations.
In fact, we have found that the performance of the full solver, in terms of total wall clock time, is actually better when a smaller time step of size 1 month is taken, due to the faster convergence of the flow problem Newton solver at every time step.

\subsection{Summary of Results}

To summarize, from the results of these three cases we conclude that the explicit fixed-stress scheme is not a reliable source of pressure stabilization in general, especially in cases representative of geologic carbon storage. A pressure stabilizing effect is only seen when large time steps are taken relative to the timescale of the solution, and this is challenging in simulations of CO$_2$ storage because (i) large time steps result in large errors in the drained regions where we expect to be performing injection, and (ii) large time steps can slow or prevent the convergence of nonlinear flow problem. Moreover, if the problem approaches a steady state, for example if injection has stopped and the CO$_2$ plume is stationary, then the explicit fixed-stress scheme cannot provide a pressure stabilizing effect for any time step size.

\section{Jump Stabilization of Pressure}

Equipped with an understanding of when the explicit fixed-stress method can fail to provide a pressure stabilizing effect, we now move on to consider the use of an additional stabilization technique in conjunction with the splitting.
There are many pressure stabilization techniques that have been developed for mixed finite element discretizations and have been used in fully implicit simulations of poromechanics \cite{truty2006stabilized, aguilar2008numerical, preisig2011stabilization, wan2003stabilized, white2008stabilized, berger2015stabilized, camargo2021macroelement,frigo2021efficient}.
A convenient choice is pressure jump stabilization \cite{hughes_1987, berger2015stabilized, berger2017stabilized, camargo2021macroelement, aronson2023pressure}, which we showed to be effective in conjunction with the iterative fixed-stress method in \cite{aronson2024FS}.

This technique is built upon the definition of an artificial flux which penalizes jumps in pressure between elements, thus smoothing the oscillatory, spurious pressure modes associated with unstable spatial discretizations.
In particular, in the compositional setting a stabilization flux for each component is defined via
\begin{equation}
    F^{n+1}_{stab, c} = \tau V \sum_{\ell=1}^{n_p} \left( x_{c\ell} \rho_\ell k_{r\ell}\right) _{upw}^{n} \left[ \left[ \frac{p^{n+1}-p^n}{\delta t}\right] \right].
    \label{eq:stab_flux}
\end{equation}
This flux augments the TPFA flux during the flow solve within each time step of the explicit fixed-stress method.
Square brackets denote the difference across a face and we use the subscript $upw$ to denote that the mass used is from the upwind cell.
The parameter $\tau$ is user-defined and can be estimated using matrix conditioning arguments as in \cite{silvester1994optimal, camargo2021macroelement, aronson2023pressure}:
\begin{equation}
    \tau = c \cdot \frac{9}{32(\lambda + 4G)}.
    \label{eq:opt_tau}
\end{equation}
Here $\lambda$ is the first Lamé parameter of the solid skeleton and $G$ is the shear modulus.
The parameter $c$ can be used to account for different element topologies.
We select $c = 1$ for hexahedra and $c = 3$ for tetrahedra, and note that these values were also shown to perform well in terms of sequential iterations required to reach convergence in the iterative fixed-stress setting in \cite{aronson2024FS}.


\subsection{Stability and Convergence for Drained Problems}

Up to this point, however, no formal analysis has been done to motivate the effectiveness of the combined pressure stabilized fixed-stress method.
In fact, interesting results for the drained setting can be obtained for simple problems using the same methodology followed in \cite{kim2011stability}.
Following that approach, we consider a one dimensional, single-phase problem discretized in space with our choice of a mixed piecewise linear - piecewise constant approach and in time with the explicit fixed-stress splitting and implicit Euler. The fully discrete flow equation for cell $j$ in this case reads
\begin{equation}
\begin{split}
  \left(\frac{1}{M \delta t} + \frac{b^2}{K_{dr} \delta t}\right) \left( p^{n+1}_j - p^n_j\right) - \left(\frac{b^2}{K_{dr} \delta t}\right) \left( p^{n}_j - p^{n-1}_j \right) + \left( \frac{b}{\delta t \delta x} \right) \left( u^n_{j+1/2} - u^n_{j-1/2} - u^{n-1}_{j+1/2} + u^{n-1}_{j-1/2} \right) \\ - \left(\frac{k}{\mu \delta x^2} \right)\left( p^{n+1}_{j+1} - 2p^{n+1}_{j} + p^{n+1}_{j-1} \right) - \left(\frac{\tau \delta x}{\delta t} \right) \left( p^{n+1}_{j+1}  - 2p^{n+1}_{j} + p^{n+1}_{j-1} - p^{n}_{j+1} + 2p^{n}_{j} - p^{n}_{j-1}  \right) = 0,
\end{split}
\label{eq:1d_discrete_mass}
\end{equation}
while the fully discrete mechanics equation can be written for node $j-1/2$ as 
\begin{equation}
    \left( -\frac{K_{dr}}{\delta x} \right) \left( u^{n+1}_{j-3/2} - 2u^{n+1}_{j-1/2} + u^{n+1}_{j+1/2} \right) - b \left( p^{n+1}_j-1 - p^{n+1}_{j} \right) = 0.
\end{equation}
The mechanics equation is identical to the one presented in \cite{kim2011stability}, as pressure jump stabilization does not modify this equation. The only difference is the inclusion of the last term within \cref{eq:1d_discrete_mass}.

We can use this formulation in a Von Neumann stability analysis, where we assume displacement and pressure solutions of the form
\begin{equation}
    p^n_j = \gamma^n e^{ij\theta}\hat{p},
\end{equation}
\begin{equation}
    u^n_j = \gamma^n e^{ij\theta}\hat{u}.
\end{equation}
Substitution of these assumed solutions yields an amplification matrix which is identical to the one derived for the fixed-stress scheme without jump stabilization \cite{kim2011stability}, except that the pressure entry of the mass equation also includes our stabilizing terms.
Solving for the amplification factors $\gamma$ which result in a singular system results in $\gamma = 0$ and 
\begin{equation}
    \gamma = \frac{\delta x^2 \mu (K_{dr}+ Mb^2) + 2 \mu K_{dr} M \delta x^3 \tau (1-\cos\theta)}{\delta x^2 \mu (K_{dr}+ Mb^2) + 2K_{dr}M \delta t k (1-\cos\theta) + 2 \mu K_{dr} M \delta x^3 \tau (1-\cos\theta)}.
\end{equation}
Note that if the stabilization parameter $\tau$ is set to zero this recovers the original amplification factor in \cite{kim2011stability}.
Interestingly, this amplification factor is actually larger than the amplification factor for the standard method, as the stabilization term $2 \mu K_{dr} M \delta x^3 \tau (1-\cos\theta)$ which is added to both the numerator and denominator is always positive for $\tau >  0$.
This is related to our choice to penalize the jumps in the pressure increment, instead of the jump in the pressure at the final time in \cref{eq:stab_flux}.
A discussion of these different options is provided in \cite{berger2015stabilized}, where it was seen that the pressure increment method was the most accurate.
If one instead penalized the pressure solution itself, the resulting amplification factor is smaller than the original found in \cite{kim2011stability}, as expected.

Regardless, the stabilized amplification factor is still always less than unity for all values of $\delta t > 0$ and $\tau > 0$, so the method with jump stabilization included remains unconditionally stable in time like the unmodified fixed-stress scheme.
Moreover, since the pressure stabilization terms would be identical in the corresponding fully implicit method, we can use the same strategy as in \cite{kim2011stability} to conclude that the pressure stabilized fixed stress scheme is also convergent with a fixed number of iterations per time step (which includes the explicit method).
In particular, the error is not magnified when compared to the pressure stabilized fully implicit method (as the error term will be identical to that presented in \cite{kim2011stability}), and it is known that the pressure stabilized fully implicit scheme is still convergent by construction.
Of course, if $\tau$ is poorly selected, and in particular is much too large, the error in the pressure stabilized method may be larger than in the original method, but it has been demonstrated that this is not the case when $\tau$ is selected according to the condition number arguments defined in \cite{camargo2021macroelement, aronson2023pressure}.

To reiterate, the above arguments regard stability and convergence of the modified splitting scheme in time in the drained regime, where there is no saddle-point structure.
We feel that these properties are useful to verify, however, so that the method can be employed in both drained and undrained regions, if desired.
We verify in the following section the robustness of the pressure stabilized fixed-stress scheme in both drained and undrained regimes.

\section{Numerical Results with Pressure Stabilization}

As a final technical section of this work, we revisit the numerical examples of Section 5, but now generate results including pressure jump stabilization.
For brevity we focus on the cantilever and HI24L examples.
We confirm that the pressure stabilization scheme has very little effect when the explicit fixed-stress solution was not oscillatory, but behaves in a similar manner to the fully implicit case when pressure oscillations were still present.   

\subsection{Cantilever Plate}

\begin{figure}
\centering
\subfloat[$t = 3$ days]{\label{sfig:cant_fim_3_stab}\includegraphics[width=.4\textwidth]{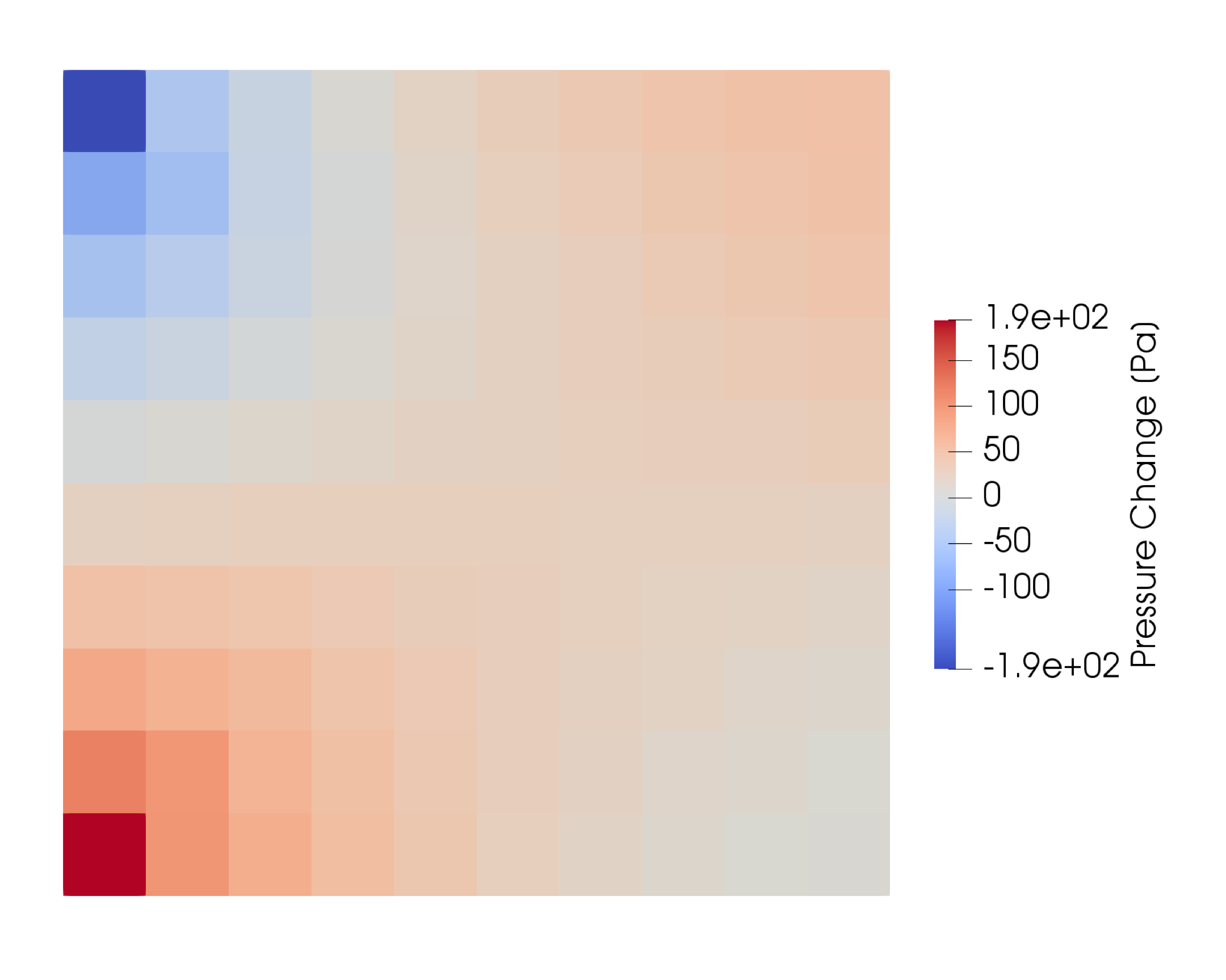}}\hfill
\subfloat[$t = 7$ days]{\label{sfig:cant_fim_7_stab}\includegraphics[width=.4\textwidth]{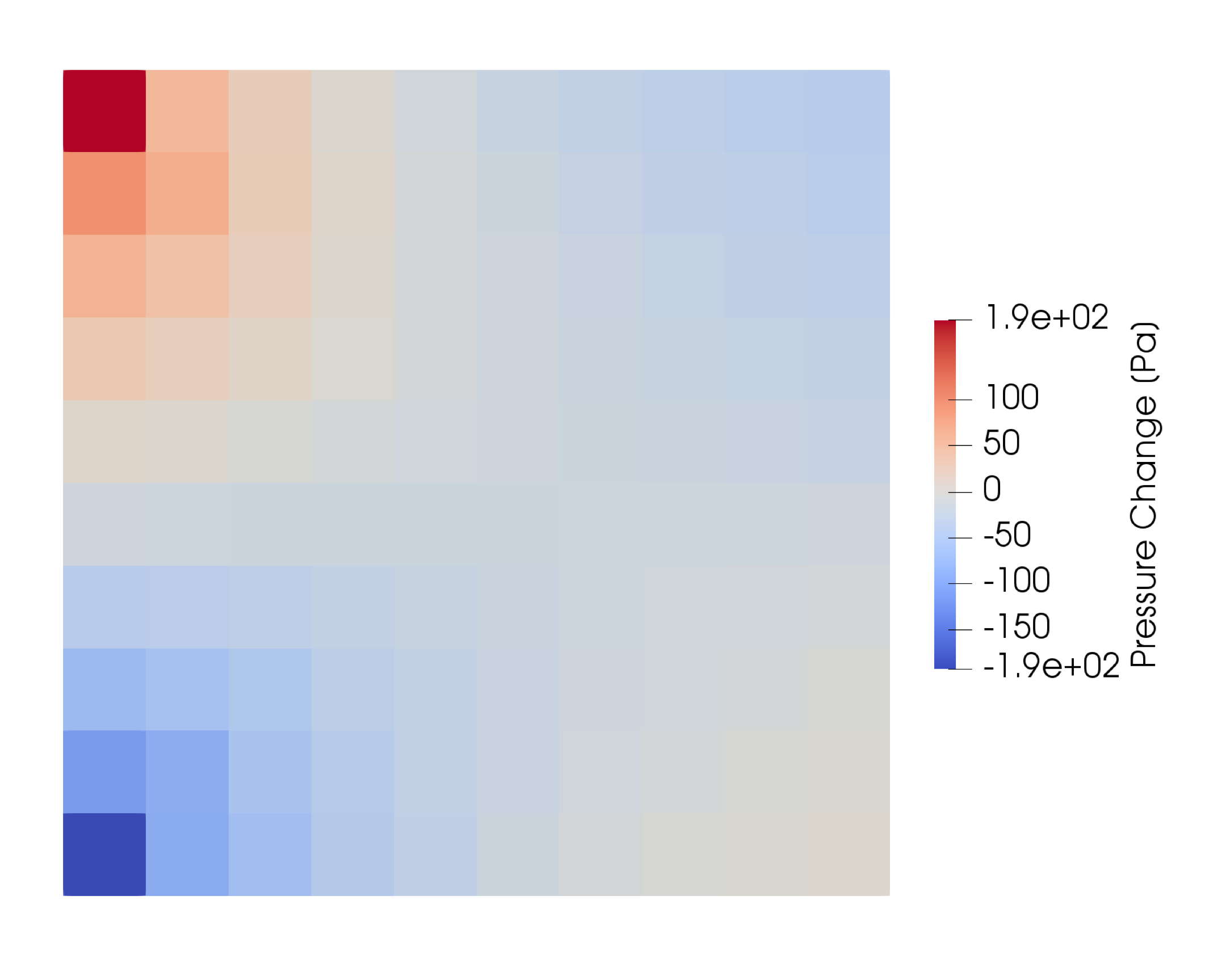}}\\
\caption{Cantilever: Pressure field at various time steps obtained with pressure stabilized fully implicit method}
\label{fig:cant_fim_stab}
\end{figure}

\begin{figure}
\centering
\subfloat[$t = 3$ days]{\label{sfig:cant_seq_1daydt_3_stab}\includegraphics[width=.4\textwidth]{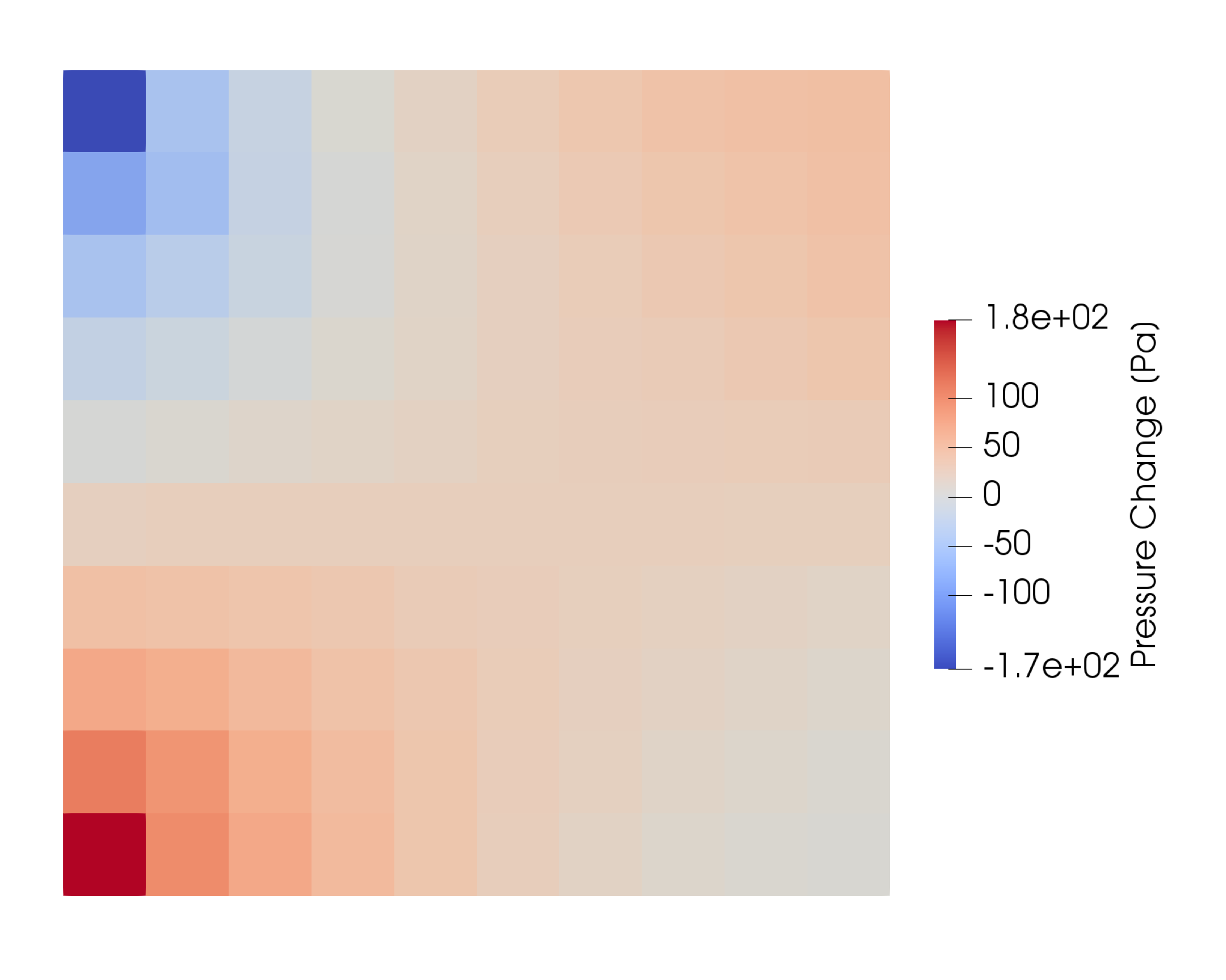}}\hfill
\subfloat[$t = 7$ days]{\label{sfig:cant_seq_1daydt_7_stab}\includegraphics[width=.4\textwidth]{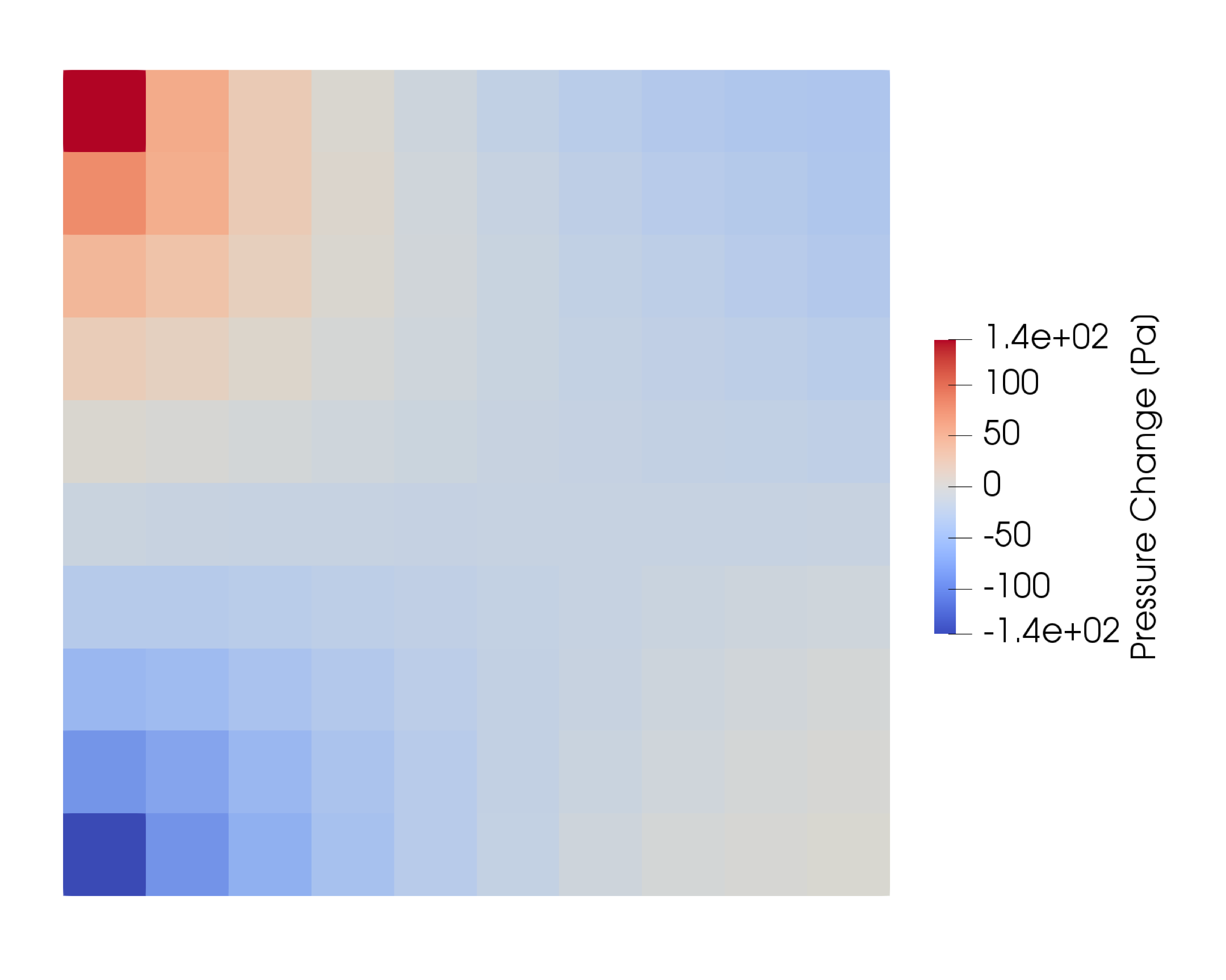}}\\
\caption{Cantilever: Pressure field at various time steps obtained with pressure stabilized explicit fixed-stress method and time step of 1 day}
\label{fig:cant_seq_1daydt_stab}
\end{figure}

\begin{figure}
\centering
\subfloat[$t = 3$ days]{\label{sfig:cant_seq_1e-2daydt_3_stab}\includegraphics[width=.4\textwidth]{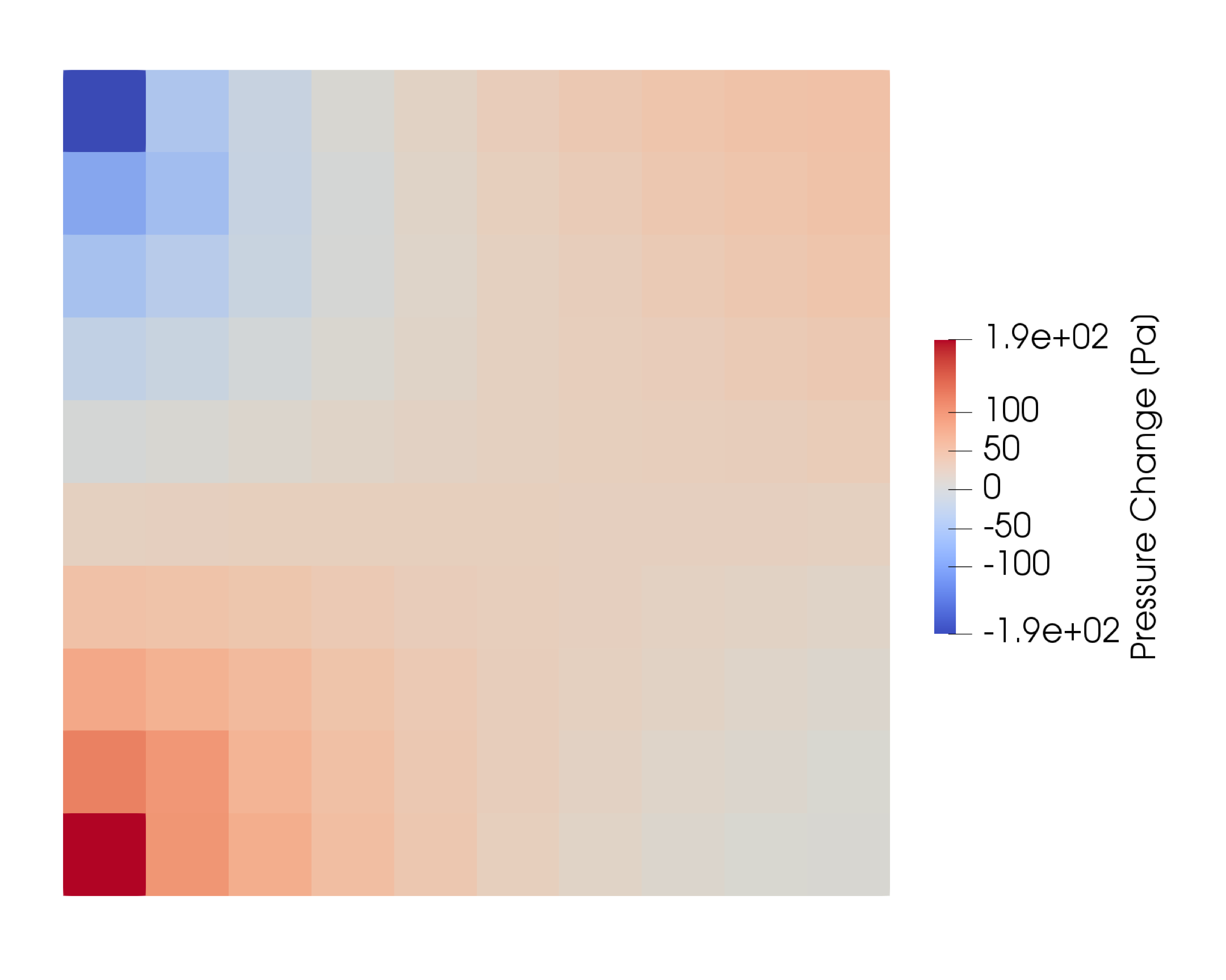}}\hfill
\subfloat[$t = 7$ days]{\label{sfig:cant_seq_1e-2daydt_7_stab}\includegraphics[width=.4\textwidth]{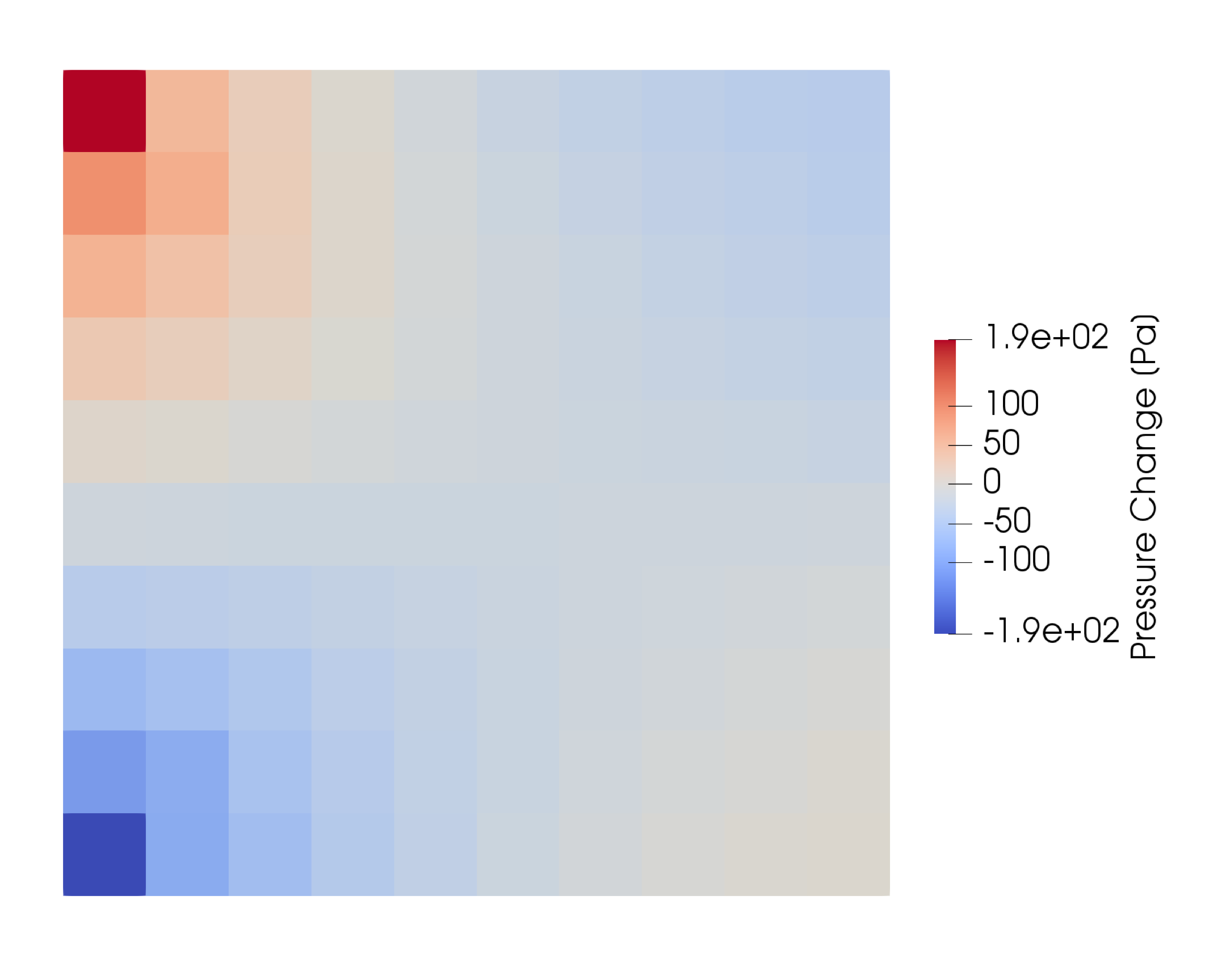}}\\
\caption{Cantilever: Pressure field at various time steps obtained with pressure stabilized explicit fixed-stress method and time step of 0.01 day}
\label{fig:cant_seq_1e-2daydt_stab}
\end{figure}

We now return to the cantilever example in which it was demonstrated that the explicit fixed-stress scheme could smooth spurious pressure oscillations in transient problems.
\cref{fig:cant_fim_stab} shows the pressure fields obtained with a fully implicit solver including jump stabilization.
Comparing this with \cref{fig:cant_seq_1daydt}, we see that the result obtained using just the explicit scheme with a time step of one day compares well with the stabilized solution.
Note that there is a difference in magnitude between the explicit and fully implicit solution, reflecting that the stabilizing effect is accompanied by error.
In \cref{fig:cant_seq_1daydt_stab} we show the results obtained with the explicit fixed-stress scheme with a time step of one day and jump stabilization included.
Indeed, the jump stabilization has affected the solution very little, meaning that it was essentially unnecessary, but it is reassuring that there was no negative effect from its inclusion.
When a time step of size 0.01 days is taken, however, \cref{fig:cant_seq_1e-2daydt_stab} illustrates that jump stabilization is as effective at smoothing the pressure oscillations as in the fully implicit case. 

\subsection{High Island 24L}

\begin{figure}
\centering
\subfloat[Fully Implicit, $c = 3$]{\label{sfig:HI_FIMC3}\includegraphics[width=.4\textwidth]{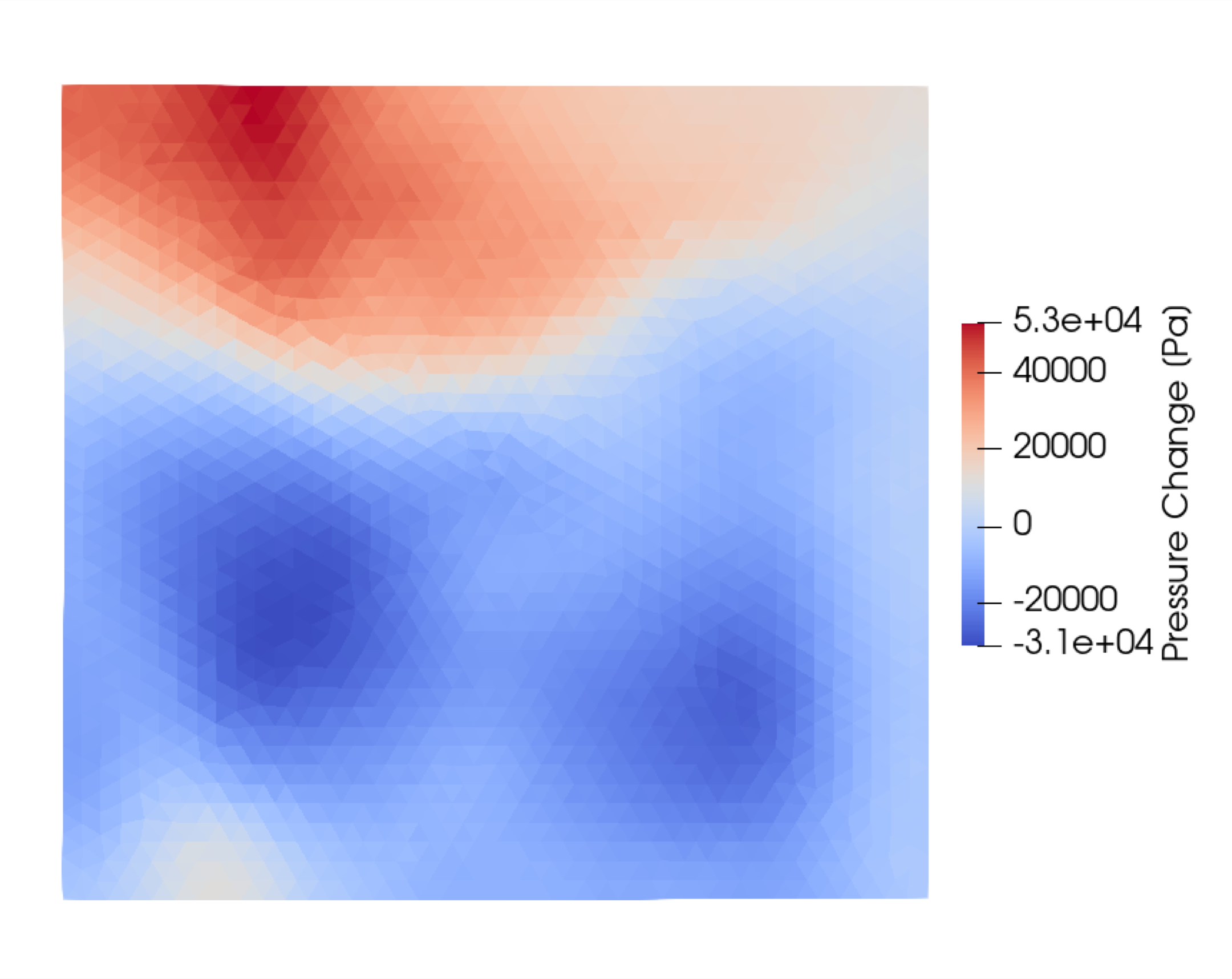}}\hfill
\subfloat[Explicit Fixed-Stress, $c = 3$]{\label{sfig:HI_SEQC3}\includegraphics[width=.4\textwidth]{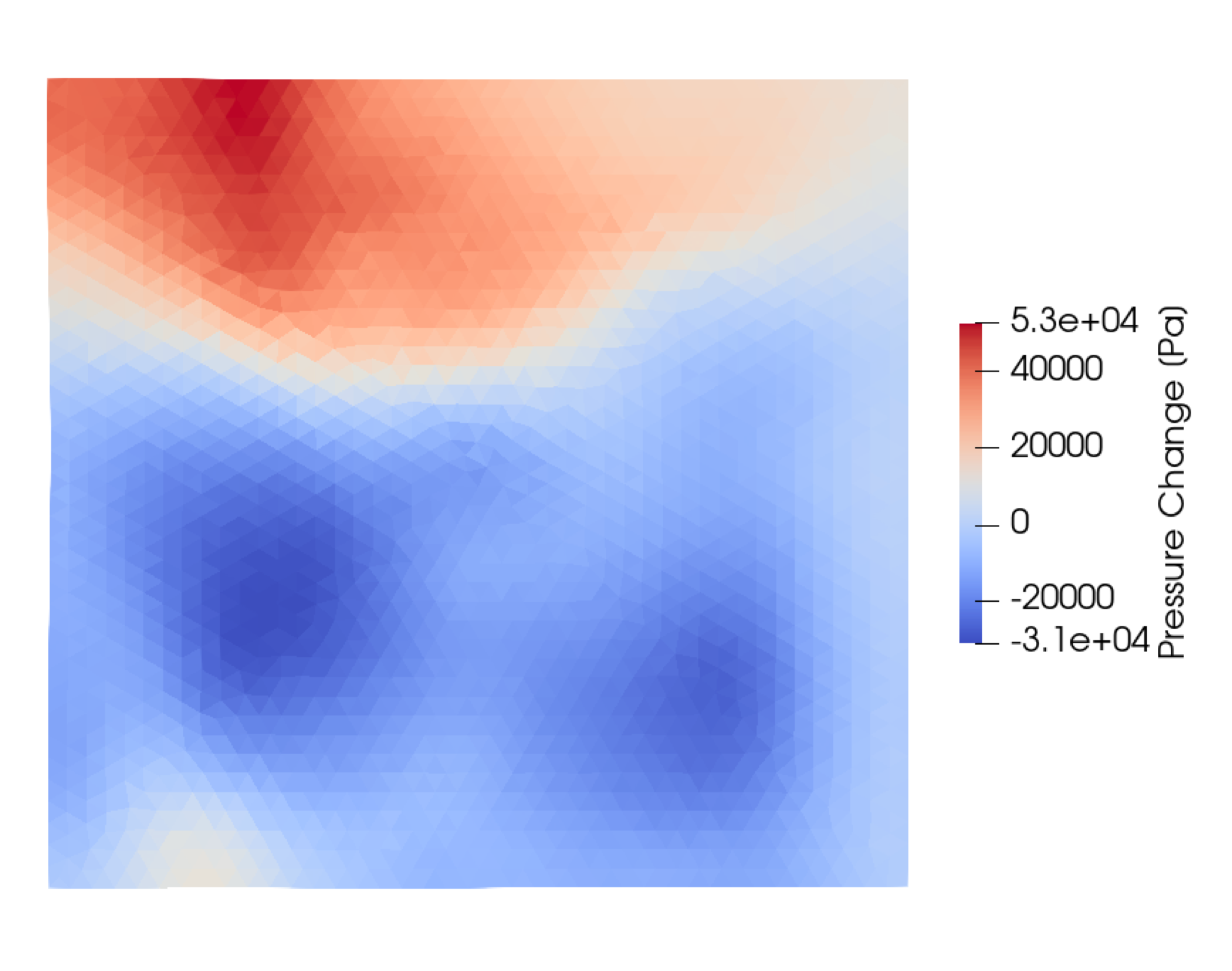}}\\
\caption{HI24L: Seafloor pressure field after 30 years of injection with time step of 6 months, including pressure stabilization}
\label{fig:HIC3}
\end{figure}

For our final set of numerical experiments, we again consider the case of CO$_2$ injection into the HI24L formation, but now with pressure jump stabilization included.
Following \cite{aronson2024FS}, we only apply stabilization in the undrained burden regions to avoid inappropriately smoothing the interface between these regions and the drained aquifer region.
\cref{fig:HIC3} shows the resulting pressure fields at the seafloor.
Clearly the spurious oscillations seen in \cref{fig:HI} have been removed, and the results obtained with the explicit fixed-stress scheme match well with the fully implicit solution.

\begin{figure}
\centering
\subfloat[Fully Implicit, $c = 0.3$]{\label{sfig:HI_FIMC03}\includegraphics[width=.4\textwidth]{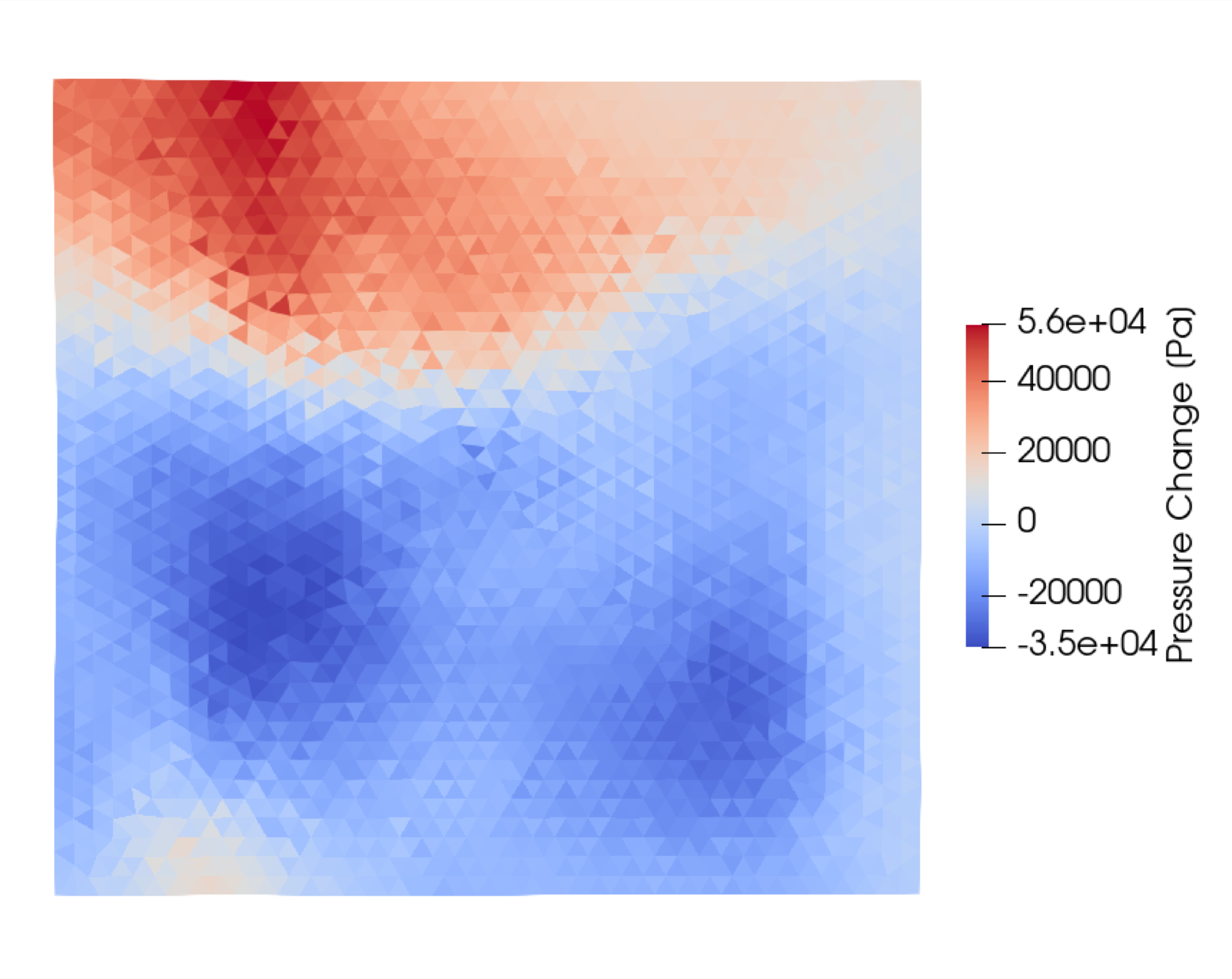}}\hfill
\subfloat[Explicit Fixed-Stress, $c = 0.3$]{\label{sfig:HI_SEQC03}\includegraphics[width=.4\textwidth]{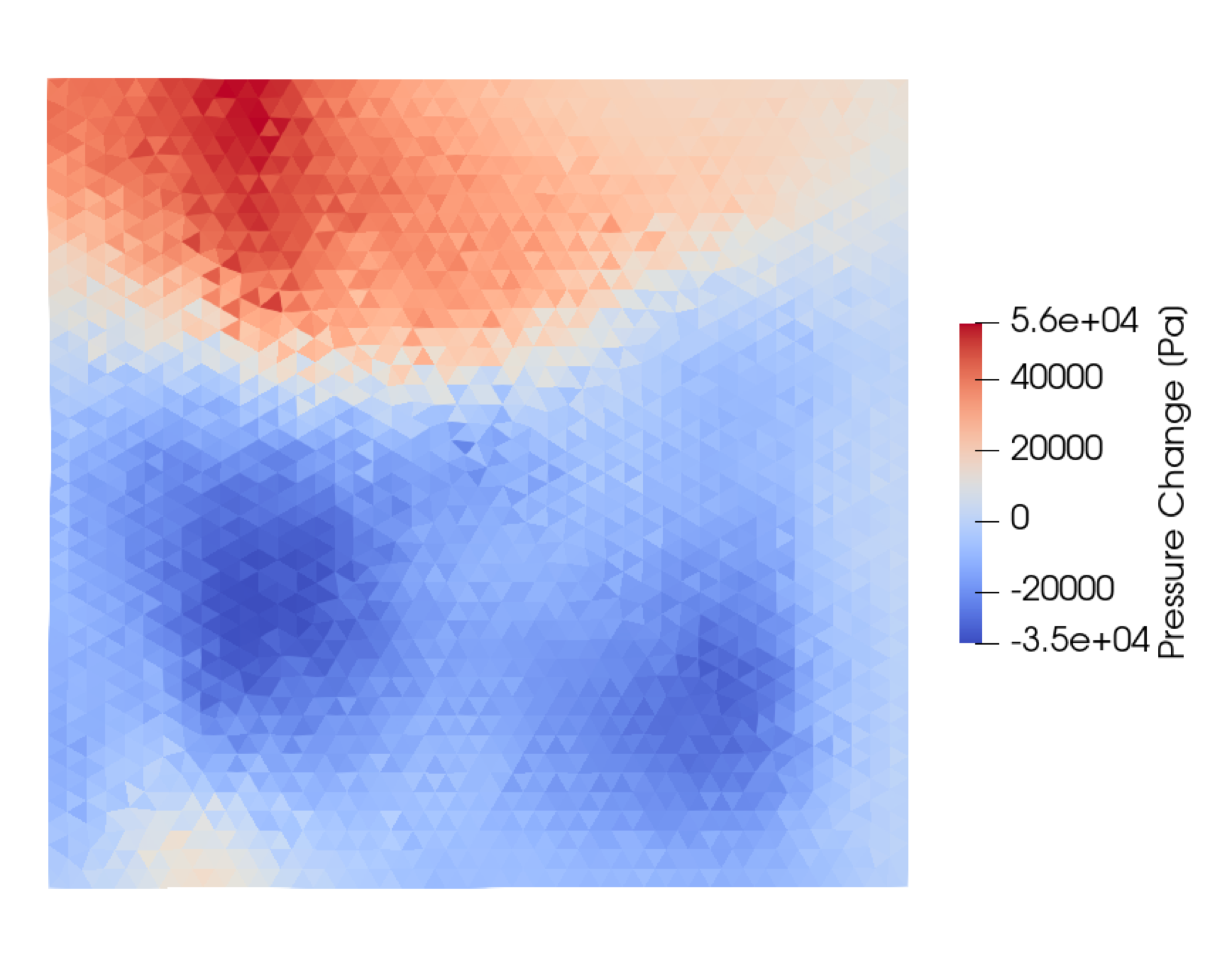}}\\
\subfloat[Fully Implicit, $c = 30$]{\label{sfig:HI_FIMC30}\includegraphics[width=.4\textwidth]{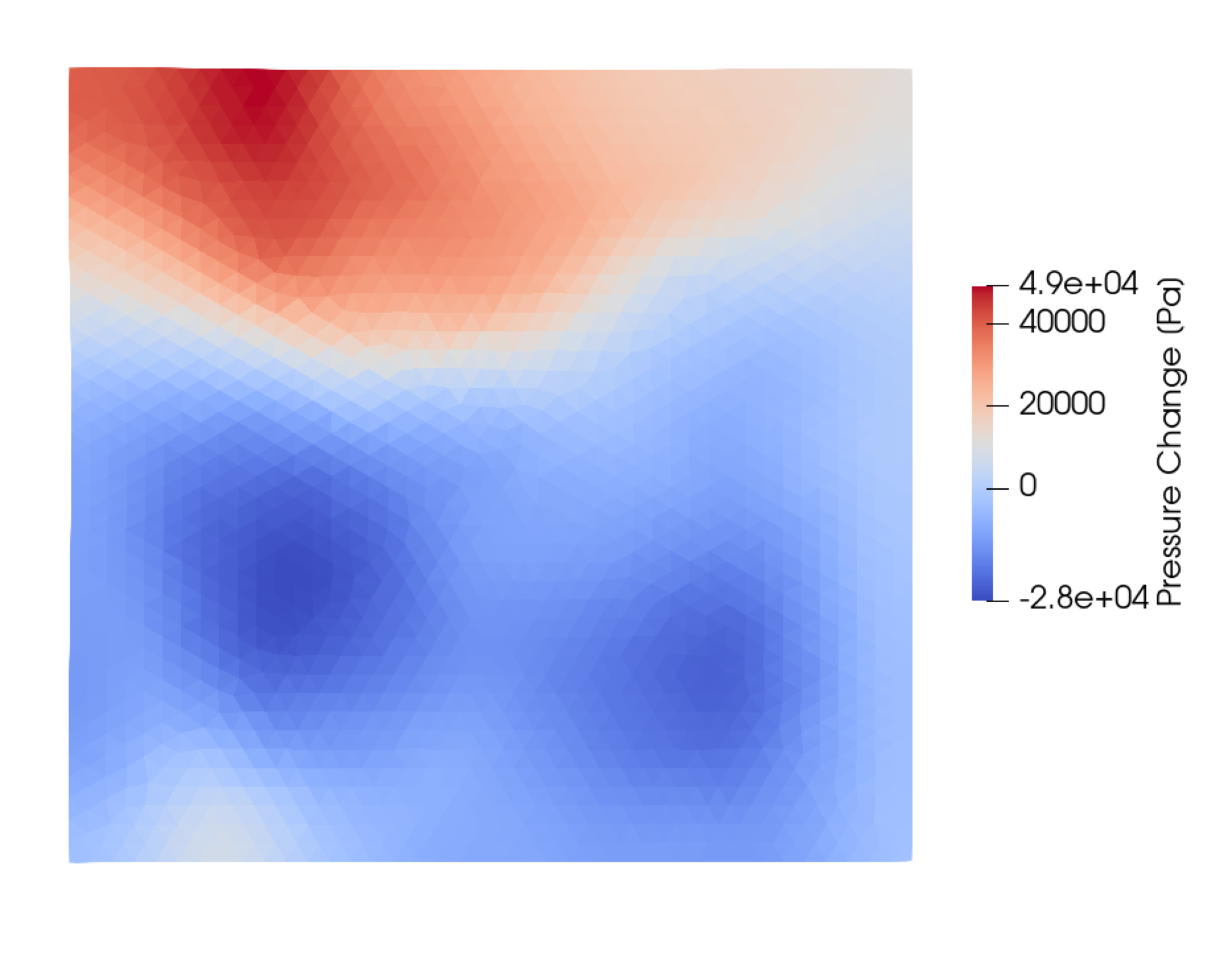}}\hfill
\subfloat[Explicit Fixed-Stress, $c = 30$]{\label{sfig:HI_SEQC30}\includegraphics[width=.4\textwidth]{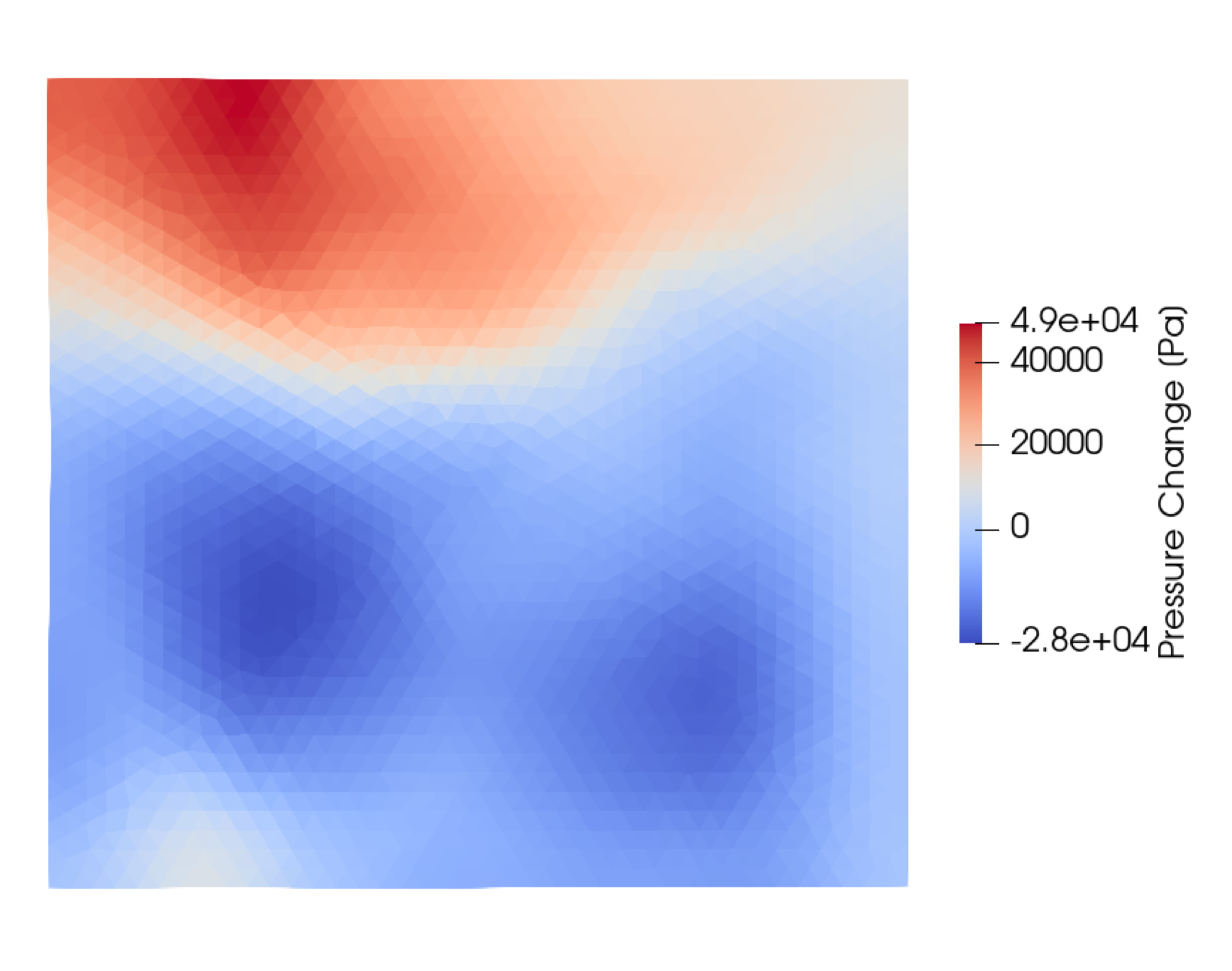}}\\
\caption{HI24L: Seafloor pressure field after 30 years of injection with time step of 6 months, including pressure stabilization of varying strengths}
\label{fig:HICvary}
\end{figure}

We also investigate the sensitivity of the results to the stabilization strength $c$ using this example.
\cref{fig:HIC3} was generated using $c = 3$, which was suggested in \cite{aronson2023pressure} to be nearly optimal in terms of conditioning of the pressure Schur complement in the fully implicit scheme, as well as in terms of iterations required for convergence of the iterative fixed-stress scheme \cite{aronson2024FS}.
For comparison, \cref{fig:HICvary} shows the seafloor pressure fields obtained with both the fully implicit method and the explicit fixed-stress, with $c = 0.3$ and $c = 30$.
We see that the sensitivity of the pressure to the parameter $c$ in the explicit fixed-stress solution is similar to what is seen in the fully implicit solution, and moreover the choice of $c = 3$ still seems to be quite good, as it smooths the spurious oscillations without over-damping the physical solution. 

\section{Conclusions}

In this work we performed detailed studies of the pressure stability of the explicit, or non-iterated, fixed-stress splitting in undrained and incompressible cases when used in conjunction with a non-inf-sup stable spatial discretization.
We detailed the similarity of the method to classical Uzawa iteration, which provided intuition on when a~pressure stabilizing effect might be expected.
The numerical studies confirmed that, for transient problems, the use of large time steps can result in pressure field visibly free of spurious oscillations.
As the time step is refined or steady state is reached, however, the spurious modes reappear.
Due to the slow evolution of the pressure field in geologic seal units, we saw that very large time steps must be taken to see the stabilizing effect.
This, however, resulted in inaccuracies in the reservoir or aquifer regions, and in the multiphase case challenged the nonlinear solver of the flow problem.
From this we concluded that, in applications such as CO$_2$ sequestration, any pressure stabilizing effect from the fixed-stress method will likely be insufficient in general.
We the discussed the use of pressure jump stabilization in conjunction with the fixed-stress method, including proofs that the addition of an additional pressure stabilization strategy does not affect the convergence or time stability of the splitting in the drained case.
By revisiting the numerical examples, we saw that the addition of jump stabilization had little effect when the oscillations were removed by the explicit fixed-stress method alone, as desired, but it also stabilized cases where the fixed-stress scheme was insufficient.
Moreover, the behavior of the scheme is the same as the pressure stabilized fully implicit method, including sensitivity to the stabilization parameter. 

\section*{Acknowledgements}
\label{sec::acknow}
Funding was provided by TotalEnergies and Chevron through the FC-MAELSTROM project.
Portions of this work were performed under the auspices of the U.S. Department of
Energy by Lawrence Livermore National Laboratory under Contract DE-AC52-07-NA27344.

\FloatBarrier

\end{document}